\documentclass[10pt]{article}

\RequirePackage{amsmath}
\RequirePackage{amsfonts}
\RequirePackage{amssymb}
\RequirePackage{hyperref}
\RequirePackage{multirow}
\RequirePackage{algorithm}
\RequirePackage{algpseudocode}
\RequirePackage{graphicx}
\RequirePackage{xcolor}
\RequirePackage{amsthm}
\RequirePackage{amsmath}
\RequirePackage{amsfonts}
\RequirePackage{amssymb}
\RequirePackage{hyperref}
\RequirePackage{multirow}
\RequirePackage{algorithm}
\RequirePackage{algpseudocode}
\RequirePackage{graphicx}
\RequirePackage{xcolor}
\allowdisplaybreaks

\usepackage{dpr_fec,ifthen,dsfont}

\usepackage{isetechreport}
\def\reportyear{26T}
\def\reportno{022}
\def\revisionno{0}
\def\originaldate{\today}
\def\revisiondate{\today}
\coraltrue
\cvcrfalse

\begin{document}

\title{A Retrospective-Approximation Interior-Point Method}

\author{Frank E.~Curtis\thanks{E-mail: frank.e.curtis@lehigh.edu}}
\author{Zahra Khatti\thanks{E-mail: zak223@lehigh.edu}}
\author{Daniel P.~Robinson\thanks{E-mail: daniel.p.robinson@lehigh.edu}}
\affil{Department of Industrial and Systems Engineering, Lehigh University}
\titlepage

\maketitle

\begin{abstract}
  An algorithm is proposed, analyzed, and tested for solving nonlinearly constrained continuous optimization problem instances wherein the objective function is defined by an expectation over a random variable argument and the constraints include inequalities (and perhaps equalities) involving functions that may be nonconvex. The algorithm combines a retrospective-approximation strategy for handling the objective and an interior-point approach for handling the constraints. In each iteration of an outer loop, the objective is approximated by an empirical average over a set of realizations of the random variable argument and the resulting sampled problem is solved with an interior-point method that allows inexact search-direction computations for computational flexibility. Termination conditions for the interior-point method are tailored for those needed for convergence guarantees by the outer loop. Over the outer loop, the sample set is updated adaptively in such a manner that convergence guarantees for solving the original optimization problem can be assured. The results of numerical experiments are presented, which show that the proposed algorithm can outperform competing methods that have comparable theoretical guarantees.
\end{abstract}


\section{Introduction}\label{sec.introduction}

The purpose of this paper is to propose, analyze, and provide the results of numerical experiments with an efficient algorithm for solving nonlinearly constrained continuous optimization problems.  In particular, our proposed algorithm is designed for settings in which the objective function is defined by an expectation over a random variable argument and the inequality (and any equality) constraints are defined by functions that may be nonconvex.  Multiple types of algorithms exist for solving such problems, each of which offers different types of convergence guarantees.  Our focus in this paper is on a class of algorithms that aims to solve the original problem to high accuracy by sampling increasingly larger numbers of terms for approximating the objective function.  At the same time, our proposed algorithm aims to exploit the fact that the original problem may be solved more efficiently by starting with a small number of sampled objective terms, employing a deterministic algorithm to solve the corresponding approximate optimization problem to moderate accuracy, and only then increasing the sample size.  This leads to an iterative scheme that aims to be efficient in terms of sampling the objective function, while also eventually ensuring that the original problem of interest is solved to high accuracy.

The algorithm we propose combines a retrospective-approximation (RA) strategy for handling the objective function and an interior-point method (IPM) for handling the constraints. RA is a paradigm for solving continuous optimization problems that has been proposed and analyzed for both unconstrained \cite{NewtBollPasuYip2025,Pasu2011} and constrained~\cite{berahas2025retrospective} settings.  It is closely related to sample-average-approximation schemes in which a sequence of sampled problems with increasing sample sizes is solved, each only approximately; see, e.g., \cite{RoysPola2004,Roys2012}, in which the sample size is increased according to rules based on optimality functions, and \cite{Roys2013}, in which the sample sizes are selected with the aim of minimizing the total computational effort.  RA schemes that employ an IPM to solve each sampled problem have been considered for two-stage stochastic linear programs in~\cite{GhosPasu2011}; our setting, with nonlinear and possibly nonconvex constraints, is different.  On the other hand, the class of IPMs is known to be one of the most effective classes of algorithmic methods for solving large-scale continuous optimization problems.  This is evidenced by IPM strategies forming the core of state-of-the-art software packages such as Ipopt \cite{wachter2006implementation}, Knitro \cite{ByrdHribNoce99}, and LOQO \cite{VandShan99}.

Our approach can be contrasted with stochastic variants of sequential quadratic programming (SQP) methods and IPMs that have been proposed for this same class of problems, i.e., with an objective defined by an expectation and nonlinear constraint functions that can be evaluated deterministically.  Stochastic SQP methods that employ a stochastic objective-gradient estimate in each iteration were first proposed for problems with nonlinear equality constraints in~\cite{BeraCurtRobiZhou2021}, and this line of work has since been extended to allow for rank-deficient constraint Jacobians through a step-decomposition strategy~\cite{BeraCurtONeiRobi2024} and to allow for inexact subproblem solutions, which, as in this paper, is of particular relevance in large-scale settings~\cite{CurtRobiZhou2024}.  For problems with nonlinear inequality and equality constraints, stochastic SQP methods include one with convergence guarantees in expectation under loose assumptions on the gradient estimates~\cite{curtis2024sequential} and an active-set method that employs an augmented Lagrangian merit function and a stochastic line search and offers almost-sure convergence guarantees~\cite{na2021inequality}.  More recently, stochastic IPMs have been proposed as well: a stochastic-gradient-based IPM for bound-constrained problems~\cite{CurtKungRobiWang2025}; a single-loop stochastic IPM for nonlinear inequality-constrained problems that maintains feasibility of the iterates with respect to the inequality constraints while the barrier parameter is driven to zero along with the stochastic-gradient iterations~\cite{curtis2026single}; and a trust-region stochastic SQP method for nonlinear inequality- and equality-constrained problems that handles the inequality constraints with an interior-point approach, decreases the barrier parameter along a prescribed sequence, and offers almost-sure convergence guarantees~\cite{FangKimNaDemmLava2026}.  Methods of this type take a stochastic step in every iteration, so their step sizes and merit or barrier parameters must be controlled to account for the stochastic error present in every iteration, and they do not directly exploit mature deterministic solvers for constrained problems.  A complementary paradigm is to increase the number of samples as the optimization proceeds, as in the progressive-batching and adaptive-sampling strategies developed for unconstrained problems~\cite{bollapragada2018progressive,bollapragada2018adaptive} and extended to equality-constrained problems in the adaptive-sampling SQP method of~\cite{BeraBollZhou2022}, so that the stochastic error is reduced as a solution is approached.  RA can be viewed as the member of this paradigm in which each sampled problem is solved by a deterministic method to a prescribed accuracy before the sample size is increased, which is precisely the feature that allows an IPM with inexact step computations to be designed.

\subsection{Contributions}

Our proposed algorithm builds on the RA framework proposed in \cite{berahas2025retrospective} for solving constrained optimization problems by combining it with an IPM that uses the inexact search direction computations proposed in \cite{CurtScheWaec2010,CurtHubeScheWaec2012}.  As opposed to the RA-SQP framework proposed in \cite{berahas2025retrospective}, which involves distinct termination criteria from the general RA framework analyzed in that paper, our RA-IPM strategy involves termination criteria for the inner IPM that are tied exactly to those for the outer RA framework; see upcoming \eqref{eq:kkt:full}, \eqref{eq:kkt:sampled}, \eqref{eq:sample:termination}, and \eqref{eq:barrier_termination}.  The results of our numerical experiments show that, amongst other advantages, our RA-IPM approach offers improved performance in terms of CPU time compared to the aforementioned RA-SQP algorithm on a set of larger-scale test problems.

\subsection{Notation}

We use $\R{}$ to denote the set of real numbers, $\R{n}$ to denote the set of $n$-dimensional real vectors, and $\R{n \times m}$ to denote the set of $n$-by-$m$-dimensional real matrices.  Given $r \in \R{}$, we use $\R{}_{>r}$ (respectively, $\R{}_{\geq r}$) to denote the set of real numbers that are greater than $r$ (respectively, greater than or equal to $r$).  We use $\N{}$ to denote the set of nonnegative integers, and for any $q \in \N{}$ with $q \geq 1$, we denote $[q] := \{1,\dots,q\}$. Given $s \in \R{q}$, its $i$th element is denoted as $s^{(i)}$ for each $i \in [q]$. Given $(a,b) \in \R{q} \times \R{q}$, we use $a \perp b$ to denote that $a_i b_i = 0$ for all $i \in [q]$.

\subsection{Organization}

Our optimization problem formulation of interest and fundamental assumptions about the problem and our proposed algorithm are provided in~\S\ref{sec.problem}.  Our proposed algorithm is presented in~\S\ref{sec.algorithm}, and our analysis of it is presented in~\S\ref{sec.analysis}.  The details of a specific implementation of our proposed algorithm are given in~\S\ref{sec.implementation}, and the results of numerical experiments that employ this implementation of our proposed algorithm are given in~\S\ref{sec.numerical}.  Concluding remarks are provided in~\S\ref{sec.conclusion}.

\section{Problem Formulations}\label{sec.problem}

Our algorithm is designed to solve nonlinearly constrained continuous optimization problems of the form
\bequation\label{eq:full}
  \min_{x \in \R{n}}\ f(x)\ \text{subject to}\ (\st) \ c_\Ecal(x) = 0\ \text{and}\ c_\Ical(x) \leq 0,
\eequation
where $f : \R{n} \to \R{}$, $c_\Ecal : \R{n} \to \R{p}$, and $c_\Ical : \R{n} \to \R{q}$ are continuously differentiable.  More precisely, it is designed for solving problems where $f$ is defined by
\bequationNN
  f(x) \equiv \E_{\xi} [F(x,\xi)],
\eequationNN
where $F : \R{n} \times \Xi \to \R{}$, $\xi$ is a random variable defined by a probability space $(\Xi, \Fcal_\xi, \P_\xi)$, and $\E_\xi$ denotes expectation taken with respect to $\P_\xi$.  Under the presumption that objective function and gradient values are not tractable to compute in practice in such a context, our proposed algorithm does not aim to solve \eqref{eq:full} directly.  Rather, as discussed shortly, it solves a sequence of approximate problems over which the objective function is approximated by an empirical average of sampled terms.

We make the following assumption about~\eqref{eq:full} throughout this work.  We note in passing that this assumption could be loosened such that the stated properties of the problem functions only need to hold almost surely over an open convex set that contains the iterates generated by any run of the algorithm. However, since our proposed algorithm and convergence guarantees follow using the same series of arguments in either case, we state the assumption as follows for the sake of brevity.

\bassumption\label{ass:full}
  The objective function $f : \R{n} \to \R{}$ is continuously differentiable and bounded below by $f_{\inf} \in \R{}$, i.e., $f(x) \geq f_{\inf}$ for all $x \in \R{n}$. In addition, its gradient function $\nabla f : \R{n} \to \R{n}$ is Lipschitz continuous and bounded, i.e., there exists $(L_{\nabla f}, \kappa_{\nabla f}) \in \R{}_{>0} \times \R{}_{>0}$ such that, for all $(x,\xbar) \in \R{n} \times \R{n}$, one has
  \bequationNN
    \|\nabla f(x) - \nabla f(\xbar)\|_2 \leq L_{\nabla f} \|x - \xbar\|_2\ \ \text{and}\ \ \|\nabla f(x)\|_2 \leq \kappa_{\nabla f}.
  \eequationNN
  The constraint functions $c_\Ecal : \R{n} \to \R{p}$ and $c_\Ical : \R{n} \to \R{q}$ are continuously differentiable and bounded, i.e., there exists $(\kappa_{c_\Ecal}, \kappa_{c_\Ical}) \in \R{}_{>0} \times \R{}_{>0}$ such that, for all $x \in \R{n}$, one has $\|c_\Ecal(x)\|_2 \leq \kappa_{c_\Ecal}$ and $\|c_\Ical(x)\|_2 \leq \kappa_{c_\Ical}$.  In addition, $\nabla c_\Ecal : \R{n} \to \R{n \times p}$ and $\nabla c_\Ical : \R{n} \to \R{n \times q}$ are Lipschitz continuous and bounded, i.e., there exist $(L_{\nabla c_\Ecal}, L_{\nabla c_\Ical}, \kappa_{\nabla c_\Ecal}, \kappa_{\nabla c_\Ical}) \in \R{}_{>0} \times \R{}_{>0} \times \R{}_{>0} \times \R{}_{>0}$ such that, for all $(x,\xbar) \in \R{n} \times \R{n}$:
  \begin{align*}
    \|\nabla c_\Ecal(x) - \nabla c_\Ecal(\xbar)\|_2 \leq L_{\nabla c_\Ecal} \|x - \xbar\|_2; &\ \ \|\nabla c_\Ecal(x)\|_2 \leq \kappa_{\nabla c_\Ecal}; \\
    \|\nabla c_\Ical(x) - \nabla c_\Ical(\xbar)\|_2 \leq L_{\nabla c_\Ical} \|x - \xbar\|_2; &\ \ \|\nabla c_\Ical(x)\|_2 \leq \kappa_{\nabla c_\Ical}.
  \end{align*}
\eassumption

Since our framework uses an IPM to solve a sequence of approximations to \eqref{eq:full}, there are benefits of using a slack reformulation so that nonnegativity only needs to be maintained for the slack variables during a run of the IPM, rather than values of $c_\Ical$ explicitly.  Thus, we introduce the following problem that is equivalent to \eqref{eq:full}:
\bequation\label{eq:full_slack}
  \min_{(x,s) \in \R{n} \times \R{q}}\ f(x)\ \st\ c_\Ecal(x) = 0,\ c_\Ical(x) + s = 0,\ \text{and}\ s \geq 0.
\eequation

Under Assumption~\ref{ass:full} and a constraint qualification, such as the Mangasarian-Fromovitz constraint qualification (MFCQ) \cite{MangFrom67}, any local minimizer of \eqref{eq:full_slack} satisfies its corresponding Karush-Kuhn-Tucker (KKT) conditions, which can be stated as
\bsubequations\label{eq:kkt:full}
  \begin{align}
    \nabla f(x) + \nabla c_\Ecal(x) y_\Ecal + \nabla c_\Ical(x) y_\Ical &= 0, \label{eq:kkt:full_dual_feas} \\
    y_\Ical - z = 0,\ \ c_\Ecal(x) = 0,\ \ c_\Ical(x) + s &= 0, \label{eq:kkt:full_primal_feas} \\
    \text{and}\ \ 0 \leq s \perp z &\geq 0. \label{eq:kkt:full_cons_comp}
  \end{align}
\esubequations
However, as is generally the case for local-search algorithms for solving nonlinearly constrained continuous optimization problems, it is also possible for the algorithm to be attracted to an infeasible point that minimizes a measure of constraint violation.  We say that a point $x$ is an infeasible stationary point with respect to \eqref{eq:full} if and only if it is infeasible for \eqref{eq:full}, yet is stationary for
\bequationNN
  \min_{x \in \R{n}}\ \thalf \|c_\Ecal(x)\|_2^2 + \thalf \|\max\{c_\Ical(x),0\}\|_2^2.
\eequationNN
Assuming that a given point $(x,s)$ has $s \geq 0$ and $c_\Ical(x) + s \geq 0$, the first-order optimality conditions with respect to this feasibility problem reduce to
\bequation\label{eq:infeasible_stationary}
  \baligned
    \nabla c_\Ecal(x) c_\Ecal(x) + \nabla c_\Ical(x) (c_\Ical(x) + s) &= 0 \ \ \text{and} \ \ 
    S (c_\Ical(x) + s) = 0,
  \ealigned
\eequation
where $S := \diag(s) \in \R{q \times q}$ has $(i,i)$-component equal to $s^{(i)}$ for all $i \in [q]$.

Our algorithm, proposed in the next section, employs a nested loop structure where the outermost loop, indexed by $k \in \N{}$, approximates the objective of \eqref{eq:full} through sampling of the random variable input $\xi$.  Specifically, each iteration $k \in \N{}$ is equipped with a prescribed sample size $N_k$ and employs a sample set
\bequation\label{eq.sample}
  \Scal_k := \{\xi_{k,1},\dots,\xi_{k,N_k}\},
\eequation
through which it approximates $f$ with $f_{\Scal_k} : \R{n} \to \R{}$ defined by
\bequationNN
  f_{\Scal_k}(x) = \frac{1}{N_k} \sum_{i \in [N_k]} F(x,\xi_{k,i}).
\eequationNN
This leads to the slack-based formulation of the sampled problem given by
\bequation\label{eq:sampled_slack}
  \min_{(x,s) \in \R{n} \times \R{q}}\ f_{\Scal_k}(x)\ \st\ c_\Ecal(x) = 0,\ c_\Ical(x) + s = 0,\ \text{and}\ s \geq 0.
\eequation
Consistent with Assumption~\ref{ass:full} for the full-sample objective function $f$, we make the following assumption throughout this work, where for simplicity we assume that the bounds and Lipschitz constant for the sampled objective function are consistent with those in Assumption~\ref{ass:full} for the full-sample objective function.

\bassumption\label{ass:sampled}
  Let $(f_{\inf}, L_{\nabla f}, \kappa_{\nabla f}) \in \R{} \times \R{}_{>0} \times \R{}_{>0}$ be defined in Assumption~\ref{ass:full}.  For all $k \in \N{}$, the sampled objective function $f_{\Scal_k} : \R{n} \to \R{}$ is continuously differentiable and bounded below by $f_{\inf}$.  In addition, for all $k \in \N{}$, the gradient function $\nabla f_{\Scal_k} : \R{n} \to \R{n}$ is Lipschitz continuous with respect to the $l_2$-norm with constant~$L_{\nabla f}$ and bounded in $l_2$-norm by $\kappa_{\nabla f}$.
\eassumption

Under Assumption~\ref{ass:sampled} and a constraint qualification such as the MFCQ, for any $k \in \N{}$ any local minimizer of \eqref{eq:sampled_slack} satisfies the corresponding KKT conditions
\bsubequations\label{eq:kkt:sampled}
  \begin{align}
    \nabla f_{\Scal_k}(x) + \nabla c_\Ecal(x) y_\Ecal + \nabla c_\Ical(x) y_\Ical &= 0, \label{eq:kkt:sampled_1} \\
    y_\Ical - z &= 0, \label{eq:kkt:sampled_2} \\
    c_\Ecal(x) &= 0, \label{eq:kkt:sampled_3} \\
    c_\Ical(x) + s &= 0, \label{eq:kkt:sampled_4} \\
    \text{and}\ \ 0 \leq s \perp z &\geq 0. \label{eq:kkt:sampled_5}
  \end{align}
\esubequations
Note that \eqref{eq:kkt:sampled_1} differs from \eqref{eq:kkt:full_dual_feas}, but \eqref{eq:kkt:sampled_2}--\eqref{eq:kkt:sampled_5} match \eqref{eq:kkt:full_primal_feas}--\eqref{eq:kkt:full_cons_comp} exactly.

\section{Algorithm}\label{sec.algorithm}

Our proposed algorithm aims to solve \eqref{eq:full} by solving a sequence of problems involving sampled objective functions.  For each sampled objective function, an IPM is employed to solve approximately the corresponding instance of~\eqref{eq:sampled_slack}. In this section, we describe the procedures involved in the outer loop as well as the IPM employed to solve each sampled problem.

\subsection{Retrospective Approximation Method}\label{sec.ra}

As mentioned, our algorithm employs a prescribed sequence of sample sizes $\{N_k\}$.  The requirements that we impose on $\{N_k\}$, as well as on other algorithmic quantities introduced in this section, are stated along with our theoretical analysis in \S\ref{sec.analysis}.

The main feature of our RA strategy is the set of termination conditions that are employed that dictate when it has been determined that a given sampled problem has been solved to sufficient accuracy and a new sample should be taken.  For each $k \in \N{}$, our algorithm employs an IPM to generate, from an initial primal-dual iterate $(x_{k,0},s_{k,0},y_{\Ecal,k,0},y_{\Ical,k,0},z_{k,0})$ and with an initial barrier parameter $\mu_{k,0}$, an approximate primal-dual solution $(x_{k,M_k},s_{k,M_k},y_{\Ecal,k,M_k},y_{\Ical,k,M_k},z_{k,M_k})$ of \eqref{eq:sampled_slack} along with a final barrier parameter $\mu_{k,M_k}$.  Here, $M_k \in \N{}$ represents the number of iterations that the IPM performs in iteration $k$; see~\S\ref{sec.interior}.  Importantly, as seen in the next section and in our theoretical analysis, our initialization procedure and termination conditions for the IPM ensure that these points satisfy $s_{k,0} > 0$, $s_{k,M_k} > 0$, $c_\Ical(x_{k,0}) + s_{k,0} \geq 0$, $c_\Ical(x_{k,M_k}) + s_{k,M_k} \geq 0$, $z_{k,0} > 0$, and $z_{k,M_k} > 0$. Consequently, determining whether an infeasible stationary point has been reached only requires checking \eqref{eq:infeasible_stationary}, and all that is required to measure the error in satisfying the KKT conditions~\eqref{eq:kkt:sampled} is to measure the error in satisfying \eqref{eq:kkt:sampled_1}--\eqref{eq:kkt:sampled_4} and the error in satisfying the orthogonality condition in \eqref{eq:kkt:sampled_5}. Following a similar, but distinct strategy from that in the RA-SQP~\cite{berahas2025retrospective}, our algorithm employs the conditions
\bsubequations\label{eq:sample:termination}
  \begin{align}
    &\ \|\nabla f_{\Scal_k}(x_{k,M_k}) + \nabla c_\Ecal(x_{k,M_k}) y_{\Ecal,k,M_k} + \nabla c_\Ical(x_{k,M_k}) y_{\Ical,k,M_k}\|_\infty \nonumber \\
    \leq&\ \gamma_{k,1} \|\nabla f_{\Scal_k}(x_{k,0}) + \nabla c_\Ecal(x_{k,0}) y_{\Ecal,k,0} + \nabla c_\Ical(x_{k,0}) y_{\Ical,k,0}\|_\infty + \epsilon_{k,1}, \label{eq:sample:termination_1}
  \end{align}
\esubequations
and
\addtocounter{equation}{-1}
\bsubequations
  \begin{align}
    \addtocounter{equation}{1}
    \|y_{\Ical,k,M_k} - z_{k,M_k}\|_\infty &\leq \gamma_{k,2} \|y_{\Ical,k,0} - z_{k,0}\|_\infty + \epsilon_{k,2}  \label{eq:sample:termination_2} \\
    \|c_\Ecal(x_{k,M_k})\|_\infty &\leq \gamma_{k,3} \|c_\Ecal(x_{k,0})\|_\infty + \epsilon_{k,3}, \label{eq:sample:termination_3} \\
    \|c_\Ical(x_{k,M_k}) + s_{k,M_k}\|_\infty &\leq \gamma_{k,4} \|c_\Ical(x_{k,0}) + s_{k,0}\|_\infty + \epsilon_{k,4}, \label{eq:sample:termination_4} \\
    \|S_{k,M_k} z_{k,M_k}\|_\infty &\leq \gamma_{k,5} \|S_{k,0} z_{k,0}\|_\infty + \epsilon_{k,5}, \label{eq:sample:termination_5}
  \end{align}
\esubequations
where $\gamma_{k,a} \in [0,1)$ and $\epsilon_{k,a} \in \R{}_{>0}$ for each $k \in \N{}$ and $a \in [5]$ are tolerances that ensure convergence in the outer loop over $k \in \N{}$; see \S\ref{sec.analysis} for further details.

We now state Algorithm~\ref{alg:ra}, which presents the outer loop of our overall algorithmic approach. The main step involves calling Algorithm~\ref{alg.ipm}, presented in \S\ref{sec.interior}.

\balgorithm[ht]
  \caption{Retrospective-Approximation Method}
  \label{alg:ra}
  \footnotesize
  \balgorithmic[1]
    \Require residual tolerances $\{\gamma_{k,a}\}_{k \in \N{}} \subset [0,1]$ and $\{\epsilon_{k,a}\}_{k \in \N{}} \subset \R{}_{>0}$ for each $a \in [5]$; increase tolerances $\{\nu_{s,k}\}_{k \in \N{}} \subset \R{}_{\geq 0}$ and $\{\nu_{z,k}\}_{k \in \N{}} \subset \R{}_{\geq 0}$; initial point $(x_0,s_0,y_{\Ecal,0},y_{\Ical,0},z_0) \in \R{n} \times \R{q} \times \R{p} \times \R{q} \times \R{q}$ satisfying $s_0 > 0$, $c_\Ical(x_0) + s_0 \geq 0$, and $z_0 > 0$; and barrier parameter $\mu_0 \in \R{}_{>0}$
    \For{$k = 0,1,2,\ldots$}
      \If{$(x_k,s_k)$ satisfies \eqref{eq:infeasible_stationary} and either $c_\Ecal(x_k) \neq 0$ or $c_\Ical(x_k) \not\leq 0$}
        \State \textbf{return} infeasible stationary point $(x_k,s_k)$
      \EndIf
      \State Sample the elements of $\Scal_k$ (see~\eqref{eq.sample}) independently from the distribution of $\xi$
      \State Set $(x_{k,0},s_{k,0},y_{\Ecal,k,0},y_{\Ical,k,0},z_{k,0},\mu_{k,0}) \gets (x_k,s_k,y_{\Ecal,k},y_{\Ical,k},z_k,\mu_k)$
      \State Call Algorithm~\ref{alg.ipm} to obtain $(x_{k,M_k},s_{k,M_k},y_{\Ecal,k,M_k},y_{\Ical,k,M_k},z_{k,M_k},\mu_{k,M_k})$ satisfying~\eqref{eq:sample:termination}
      \State Choose $\iota_{s,k} \in \R{q}_{\geq 0}$ with $\|\iota_{s,k}\|_\infty \leq \nu_{s,k}$, $\iota_{z,k} \in \R{q}_{\geq 0}$ with $\|\iota_{z,k}\|_\infty \leq \nu_{z,k}$, and $\iota_{\mu,k} \in [0,\infty)$
      \State Set \label{step.slack_perturbation}
      \begin{multline*}
        (x_{k+1},s_{k+1},y_{\Ecal,k+1},y_{\Ical,k+1},z_{k+1},\mu_{k+1}) \\ \gets (x_{k,M_k}, s_{k,M_k} + \iota_{s,k}, y_{\Ecal,k,M_k}, y_{\Ical,k,M_k}, z_{k,M_k} + \iota_{z,k}, \mu_{k,M_k} + \iota_{\mu,k})
      \end{multline*}
    \EndFor
  \ealgorithmic
\ealgorithm

Our IPM in the next subsection is guaranteed, under relatively loose assumptions, to ensure that either an iterate satisfying the termination conditions~\eqref{eq:sample:termination} is obtained or that, in the limit, feasibility with respect to \eqref{eq:sampled_slack} is not attained while the error in satisfying the stationarity conditions~\eqref{eq:infeasible_stationary} vanishes.  This latter possibility is handled by inner termination conditions that we present in~\S\ref{sec.interior}.

\subsection{Interior-Point Method (IPM)}\label{sec.interior}

The IPM that we propose to employ for solving (approximately) each instance of \eqref{eq:sampled_slack} is essentially that proposed in \cite{CurtScheWaec2010}.  (We also direct the reader to \cite{CurtHubeScheWaec2012}, which presents additional implementation features that can improve the practical performance of the algorithm.  However, in our discussion here we discuss the method from \cite{CurtScheWaec2010}, since it is more straightforward to present.) For completeness, we present the complete details of the method, although the main distinguishing feature of the method that we propose here and that in \cite{CurtScheWaec2010} is the set of termination conditions that we employ.  In particular, we introduce termination conditions that are tailored for the outer loop stated in Algorithm~\ref{alg:ra}.

As is typical for an IPM, our algorithm for solving \eqref{eq:sampled_slack} employs a nested loop where an outer loop updates a barrier parameter and an inner loop aims at solving a corresponding barrier subproblem.  Letting $\mu_{k,0} \in \R{}_{>0}$ denote an initial value of the barrier parameter corresponding to outer loop iteration $k$ (i.e., for solving the $k$th sampled problem~\eqref{eq:sampled_slack}), a loop for updating the barrier parameter, indexed by $j \in \N{}$, has in iteration $j$ a value $\mu_{k,j} \in \R{}_{>0}$ of the barrier parameter.  An inner loop, indexed by $l \in \N{}$, aims to solve the barrier subproblem
\bequation\label{eq:barrier_subproblem}
  \baligned
    \min_{(x,s) \in \R{n} \times \R{q}} &\ f_{\Scal_k}(x) - \mu_{k,j} \sum_{i \in [q]} \ln (s^{(i)}) \\
    \st &\ c_\Ecal(x) = 0,\ c_\Ical(x) + s = 0,\ \text{and}\ s > 0.
  \ealigned
\eequation
(Here, the constraint $s > 0$ is included to emphasize that the slack variables need to be positive so that the logarithms evaluate to real numbers, and that the algorithm needs to maintain positivity of the slack variables.)  In iteration $l$, the algorithm has a primal-dual iterate given by $(x_{k,j,l},s_{k,j,l},y_{\Ecal,k,j,l},y_{\Ical,k,j,l},z_{k,j,l})$. Let us now describe the search direction computation and iterate update procedures employed in the inner loop over the index $l$.  If these procedures lead to an iterate that is an infeasible stationary point, then the IPM terminates.  Otherwise, once \eqref{eq:barrier_subproblem} has been solved to sufficient accuracy, the method either invokes a decrease in the barrier parameter or termination of outer iteration $k$ in Algorithm~\ref{alg:ra}.

The search direction computation in iteration $(k,j,l)$ employs scaled derivatives, as described in \cite{CurtScheWaec2010}.  Specifically, let us define the scaled first-order derivative values
\bequationNN
  g_{k,j,l} := \bbmatrix \nabla f_{\Scal_k}(x_{k,j,l}) \\ -\mu_{k,j} e \ebmatrix \in\R{n+q} \ \ \text{and}\ \ J_{k,j,l} := \bbmatrix \nabla c_\Ecal(x_{k,j,l})^T & 0 \\ \nabla c_\Ical(x_{k,j,l})^T & S_{k,j,l} \ebmatrix \in\R{(p+q)\times(n+q)}
\eequationNN
as well as the symmetric matrix
\bequationNN
  W_{k,j,l} := \bbmatrix H_{k,j,l} & 0 \\ 0 & S_{k,j,l} Z_{k,j,l} \ebmatrix \in \R{(n+q)\times(n+q)}\ \ \text{where}\ \ H_{k,j,l} \in \R{n\times n}
\eequationNN
and $Z_{k,j,l} = \diag(z_{k,j,l})$.  Here, to obtain fast local convergence near a local minimizer of \eqref{eq:sampled_slack}, the symmetric matrix $H_{k,j,l}$ could be set as the Hessian of the Lagrangian of \eqref{eq:sampled_slack} with respect to the variable $x$ at the current primal-dual iterate, or an approximation of this Hessian.  Moreover, in place of $S_{k,j,l} Z_{k,j,l}$, one might use  $\mu_{k,j} I$; see \cite{CurtScheWaec2010}.  For our purposes, since we do not focus on fast local convergence guarantees, we merely presume that a Hessian modification strategy is employed, as in \cite{CurtScheWaec2010}.

Iteration $(k,j,l)$ commences by computing a normal step toward satisfaction of a linearization of the constraints at the current iterate.  Letting
\bequationNN
  c_{k,j,l} := \bbmatrix c_\Ecal(x_{k,j,l}) \\ c_\Ical(x_{k,j,l}) + s_{k,j,l} \ebmatrix,
\eequationNN
the normal step subproblem is defined as
\bequation\label{eq:normal_step}
  \min_{v \in \R{n+q}}\ \thalf \|c_{k,j,l} + J_{k,j,l} v\|_2^2\ \st\ \|v\|_2 \leq \omega \|J_{k,j,l}^T c_{k,j,l}\|_2,
\eequation
where $\omega \in \R{}_{>0}$ is a user-defined constant.  The vector $v_{k,j,l} \in \R{n + q}$ computed by the algorithm does not need to solve \eqref{eq:normal_step} exactly.  Rather, it only needs to satisfy the following loose condition, which is common in continuous optimization methods.

\smallskip
\noindent \textbf{Normal Step Condition} \textit{The vector $v_{k,j,l}$ is feasible for \eqref{eq:normal_step} and, for some user-defined $\epsilon_v \in (0,1]$, satisfies the Cauchy decrease condition
\bequationNN
  \|c_{k,j,l}\|_2 - \|c_{k,j,l} + J_{k,j,l} v_{k,j,l}\|_2 \geq \epsilon_v (\|c_{k,j,l}\|_2 - \|c_{k,j,l} + J_{k,j,l} v_{k,j,l}^C\|_2),
\eequationNN
where $v_{k,j,l}^C$ is the Cauchy point with respect to problem~\eqref{eq:normal_step}; see \cite{CurtScheWaec2010}.}
\smallskip

The algorithm next computes a tangential step $u_{k,j,l}$, which is done through computation of the combined step $d_{k,j,l} = v_{k,j,l} + u_{k,j,l}$.  In the extreme case, with
\bequationNN
  y_{k,j,l} := \bbmatrix y_{\Ecal,k,j,l} \\ y_{\Ical,k,j,l} \ebmatrix,
\eequationNN
the algorithm computes this step by solving the linear system of equations
\bequation\label{eq.linsys}
  \bbmatrix W_{k,j,l} & J_{k,j,l}^T \\ J_{k,j,l} & 0 \ebmatrix \bbmatrix d_{k,j,l} \\ \delta_{y,k,j,l} \ebmatrix = - \bbmatrix g_{k,j,l} + J_{k,j,l}^Ty_{k,j,l} \\ -J_{k,j,l} v_{k,j,l} \ebmatrix.
\eequation
However, to allow flexibility, the algorithm merely requires $(d_{k,j,l},\delta_{y,k,j,l})$ to be computed such that, with residuals $\rho_{k,j,l} \in \R{n + q}$ and $r_{k,j,l} \in \R{p + q}$ defined by
\bequation\label{eq.residuals}
  \bbmatrix \rho_{k,j,l} \\ r_{k,j,l} \ebmatrix := \bbmatrix W_{k,j,l} & J_{k,j,l}^T \\ J_{k,j,l} & 0 \ebmatrix \bbmatrix d_{k,j,l} \\ \delta_{y,k,j,l} \ebmatrix + \bbmatrix g_{k,j,l} + J_{k,j,l}^Ty_{k,j,l} \\ -J_{k,j,l} v_{k,j,l} \ebmatrix
\eequation
one of three termination tests is satisfied. These termination tests are stated in terms of three intermediate conditions.  The first relates to the first residual component.

\smallskip
\noindent \textbf{Residual Condition} \textit{The residual $\rho_{k,j,l}$ must satisfy, for some $\kappa \in (0,1)$,
\bequationNN
  \|\rho_{k,j,l}\|_2 \leq \kappa \min \left\{ \left\| \bbmatrix g_{k,j,l} + J_{k,j,l}^Ty_{k,j,l} \\ -J_{k,j,l} v_{k,j,l} \ebmatrix \right\|_2, \left\| \bbmatrix g_{k,j,l-1} + J_{k,j,l-1}^Ty_{k,j,l} \\ -J_{k,j,l-1} v_{k,j,l-1} \ebmatrix \right\|_2 \right\}.
\eequationNN}
\smallskip

In this Residual Condition, the first term in the minimum on the right-hand side is typical for inexact Newton-type methods. The second term is included to ensure convergence, as it ensures that in certain conditions a contraction occurs; see \cite{CurtScheWaec2010}.

The second intermediate condition requires either that the tangential step is small in norm relative to the normal step or that the curvature along the tangential step is sufficiently large and the objective of a tangential step subproblem (see \cite[eq.~(2.12)]{CurtScheWaec2010}) is sufficiently small relative to the norm of the normal step. This latter condition, required in the convergence analysis, ensures that any increase in the subproblem objective function is at most proportional to the norm of the normal step.

\smallskip
\noindent \textbf{Tangential Step Condition} \textit{The tangential step $u_{k,j,l} = d_{k,j,l} - v_{k,j,l}$ must satisfy, for some user-defined constant $\psi \in \R{}_{>0}$, either the inequality $\|u_{k,j,l}\|_2 \leq \psi \|v_{k,j,l}\|_2$ or, for some user-defined constants $\theta \in \R{}_{>0}$ and $\zeta \in \R{}_{>0}$, the inequalities
\begin{align*}
  \thalf u_{k,j,l}^T W_{k,j,l} u_{k,j,l} &\geq \theta \|u_{k,j,l}\|_2^2 \\
  \text{and}\ \ (g_{k,j,l} + W_{k,j,l} v_{k,j,l})^T u_{k,j,l} + \thalf u_{k,j,l}^T W_{k,j,l} u_{k,j,l} &\leq \zeta \|v_{k,j,l}\|_2.
\end{align*}}
\smallskip

In order to ensure that the Tangential Step Condition may ultimately be satisfied, the algorithm employs the following Hessian modification strategy; see \cite{CurtScheWaec2010}.

\smallskip
\noindent \textbf{Hessian Modification Strategy} \textit{Let $W_{k,j,l}$ be the current value of this matrix employed in \eqref{eq.linsys}--\eqref{eq.residuals} and let $(d_{k,j,l},\delta_{y,k,j,l})$ be a given trial step.  If $\|u_{k,j,l}\|_2 \leq \psi \|v_{k,j,l}\|_2$ or $\thalf u_{k,j,l}^T W_{k,j,l} u_{k,j,l} \geq \theta \|u_{k,j,l}\|_2^2$, where $(\psi,\theta)$ are defined within the Tangential Step Condition, then maintain this value for the matrix; otherwise, modify $W_{k,j,l}$ to increase its smallest eigenvalue by at least a positive constant amount.}
\smallskip

The third intermediate condition relates to a reduction in a local model of a merit function.  For these purposes, with respect to any barrier parameter $\mu_{k,j}$ and a positive merit parameter $\pi$, the value of which is updated adaptively by the algorithm, we introduce the merit function $\phi_{\Scal_k}(\cdot,\cdot;\mu_{k,j},\pi) : \R{n} \times \R{q} \to \R{}$ defined by
\bequationNN
  \phi_{\Scal_k}(x,s ; \mu_{k,j},\pi) = f_{\Scal_k}(x) - \mu_{k,j} \sum_{i \in [q]} \ln(s^{(i)}) + \pi \left\| \bbmatrix c_\Ecal(x) \\ c_\Ical(x) + s \ebmatrix \right\|_2.
\eequationNN
Given $(\mu_{k,j},\pi)$, a local model of this merit function at a given iterate $(x_{k,j,l},s_{k,j,l})$ is given by $m_{k,j,l}(\cdot; \mu_{k,j},\pi) : \R{n+q} \to \R{}$ as defined by
\bequationNN
  m_{k,j,l}(d;\mu_{k,j},\pi) = f_{\Scal_k}(x_{k,j,l}) - \mu_{k,j} \sum_{i \in [q]} \ln(s_{k,j,l}^{(i)}) + g_{k,j,l}^Td + \pi \| c_{k,j,l} + J_{k,j,l} d \|_2,
\eequationNN
meaning that the reduction in this local model for a given $d_{k,j,l}$ can be measured by
\begin{align*}
  \Delta m_{k,j,l}(d_{k,j,l}; \mu_{k,j},\pi)
    :=&\ m_{k,j,l}(0; \mu_{k,j},\pi) - m_{k,j,l}(d_{k,j,l}; \mu_{k,j},\pi) \\
    =&\ -g_{k,j,l}^Td_{k,j,l} + \pi (\|c_{k,j,l}\|_2 - \| c_{k,j,l} + J_{k,j,l} d_{k,j,l} \|_2).
\end{align*}
The third intermediate condition can now be stated as follows, which refers to the merit parameter values $\pi_{k,j,l-1}$ and $\pi_{k,j,l}$ that are computed within the algorithm.

\smallskip
\noindent \textbf{Model Reduction Condition} \textit{For some $\pi_{k,j,l} \geq \pi_{k,j,l-1} > 0$, the full step $d_{k,j,l}$ must satisfy, for some user-defined constant $\sigma \in (0,1)$ and $\theta$ defined within the Tangential Step Condition, the inequality
\begin{align*}
  \Delta m_{k,j,l}(d_{k,j,l}; \mu_{k,j},\pi_{k,j,l}) \geq&\ \max\{ \thalf u_{k,j,l}^T W_{k,j,l} u_{k,j,l}, \theta \|u_{k,j,l}\|_2^2 \} \\
  &\ + \sigma \pi_{k,j,l} (\|c_{k,j,l}\|_2 - \| c_{k,j,l} + J_{k,j,l} v_{k,j,l} \|_2).
\end{align*}}

We are now prepared to state the three termination tests that the algorithm employs for computing a full step $d_{k,j,l} = v_{k,j,l} + u_{k,j,l}$ composed of the normal step $v_{k,j,l}$ and a tangential step $u_{k,j,l}$ as well as a dual step $\delta_{y,k,j,l}$.  Each test states a set of conditions that determine a step to be acceptable.  The IPM thus employs an iterative method for solving~\eqref{eq.linsys}, terminating as soon as one of the termination tests states that an acceptable primal-dual step has been obtained.

The first termination test essentially ensures that a step has been computed that corresponds to a sufficient reduction in the local model of the merit function for the most recent value of the merit parameter.

\smallskip
\noindent \textbf{Termination Test 1 (TT1)} \textit{A step $(d_{k,j,l},\delta_{y,k,j,l})$ is acceptable if the Residual Condition holds, the Tangential Step Condition holds, and the Model Reduction Condition holds with $\pi_{k,j,l} = \pi_{k,j,l-1}$.}
\smallskip

The second termination test is designed for situations in which the primal iterate is close to stationarity, and an update in the dual space is sufficiently productive.

\smallskip
\noindent \textbf{Termination Test 2 (TT2)} \textit{If, for some user-defined constant $\epsilon_2 \in \R{}_{>0}$, one finds
\bequationNN
  \|J_{k,j,l}^T c_{k,j,l}\|_2 \leq \epsilon_2 \| g_{k,j,l} + J_{k,j,l}^T y_{k,j,l} \|_2,
\eequationNN
then a step $(d_{k,j,l},\delta_{y,k,j,l}) \gets (0,\delta_{y,k,j,l})$ is acceptable if, with the normal step reset to zero, i.e., $v_{k,j,l} \gets 0$, the Residual Condition holds with $d_{k,j,l} \gets 0$.}
\smallskip

The third termination test is designed for situations in which further progress in the primal space requires an increase in the merit parameter.

\smallskip
\noindent \textbf{Termination Test 3 (TT3)} \textit{A step $(d_{k,j,l},\delta_{y,k,j,l})$ is acceptable if the Residual Condition holds, the Tangential Step Condition holds, and, for some $\epsilon_3 \in (0,1)$,
\bequationNN
  \|c_{k,j,l}\|_2 - \|c_{k,j,l} + J_{k,j,l} d_{k,j,l}\|_2
  \geq
  \epsilon_3(\|c_{k,j,l}\|_2 - \|c_{k,j,l} + J_{k,j,l} v_{k,j,l}\|_2)
  > 0.
\eequationNN}
\smallskip

In case an acceptable step is computed satisfying either Termination Test 1 or~2, then the algorithm sets the new merit parameter value as $\pi_{k,j,l} \gets \pi_{k,j,l-1}$.  Otherwise, the new merit parameter value should be set to ensure that the computed search direction is one of sufficient descent for the merit function.  This is accomplished by setting, for some user-defined constant $\tau \in (0,1)$, the trial value
\bequation\label{eq.pi_trial}
  \hat\pi_{k,j,l} \gets \frac{g_{k,j,l}^Td_{k,j,l} + \max\{\thalf u_{k,j,l}^TW_{k,j,l} u_{k,j,l}, \theta \|u_{k,j,l}\|_2^2\}}{(1 - \tau) (\|c_{k,j,l}\|_2 - \|c_{k,j,l} + J_{k,j,l} d_{k,j,l}\|_2)}
\eequation
then setting the new merit parameter $\pi_{k,j,l}$ as equal to the maximum of $\pi_{k,j,l-1}$ and this trial merit parameter value plus a user-defined constant $\epsilon_\pi \in \R{}_{>0}$.

All that remains to describe iteration $(k,j,l)$ of our IPM is to describe the procedures for updating the primal-dual iterate and for testing for termination.  Due to the algorithm's use of scaled derivatives, the primal step is
\bequationNN
  \bbmatrix \dtilde_{x,k,j,l} \\ \dtilde_{s,k,j,l} \ebmatrix := \dtilde_{k,j,l} := \bbmatrix d_{x,k,j,l} \\ S_{k,j,l} d_{s,k,j,l} \ebmatrix,\ \text{where}\ \bbmatrix d_{x,k,j,l} \\ d_{s,k,j,l} \ebmatrix := d_{k,j,l}
\eequationNN
decomposes the computed search direction $d_{k,j,l}$ into the components corresponding to the $x$ and $s$ variables.  Thus, the step sizes to be computed are also defined with respect to the search direction components $\dtilde_{x,k,j,l} \in \R{n}$ and $\dtilde_{s,k,j,l} \in \R{q}$.

First, as is standard for an IPM, the algorithm uses a fraction-to-the-boundary rule to maintain positivity of the slack variables.  This is enforced by computing, for some user-defined constant $\eta_1 \in (0,1)$, the step-size limit
\bequation\label{eq:alpha_max}
  \alpha_{k,j,l}^{\max} = \min \left\{ 1, \min_{\left\{i: \dtilde_{s,k,j,l}^{(i)} < 0 \right\}} \frac{-\eta_1 s^{(i)}_{k,j,l}}{\dtilde_{s,k,j,l}^{(i)}} \right\},
\eequation
from which it follows that $\alpha_{k,j,l} \in (0,\alpha_{k,j,l}^{\max}]$ will have
\bequation\label{eq:fraction_to_boundary}
  s_{k,j,l} + \alpha_{k,j,l} \dtilde_{s,k,j,l} \geq (1-\eta_1) s_{k,j,l}.
\eequation
Second, the algorithm performs a backtracking line search to ensure that, for some user-defined constant $\eta_2 \in (0,1)$, the following sufficient decrease condition holds:
\bequation\label{eq:armijo}
  \baligned
    &\ \phi_{\Scal_k}(x_{k,j,l} + \alpha_{k,j,l} \dtilde_{x,k,j,l}, s_{k,j,l} + \alpha_{k,j,l} \dtilde_{s,k,j,l}; \mu_{k,j}, \pi_{k,j,l}) \\
    \leq &\ \phi_{\Scal_k}(x_{k,j,l}, s_{k,j,l}; \mu_{k,j}, \pi_{k,j,l}) - \eta_2 \alpha_{k,j,l} \Delta m_{k,j,l}(d_{k,j,l}; \mu_{k,j}, \pi_{k,j,l}).
  \ealigned
\eequation
Third, to ensure that the iterations maintain $c_\Ical(x_{k,j,l}) + s_{k,j,l} \geq 0$ for all $(k,j,l)$, the algorithm performs a slack reset, which involves the component-wise update
\bequation\label{eq:slack_reset}
  s_{k,j,l+1} \gets \max\{s_{k,j,l+1},-c_\Ical(x_{k,j,l+1})\}.
\eequation

As for the update in the dual space, the algorithm needs to update the multipliers for the nonlinear constraints as well as the multipliers for the bounds on the slack variables.  With respect to the former, the algorithm merely needs to ensure that a new stacked multiplier vector $y_{k,j,l+1}$ is computed so as to satisfy the inequality
\bequation\label{eq.mystery_inequality}
  \|g_{k,j,l} + J_{k,j,l}^T y_{k,j,l+1}\|_2 \leq \|g_{k,j,l} + J_{k,j,l}^T(y_{k,j,l} + \delta_{y,k,j,l})\|_2.
\eequation
To be more specific for the sake of clarity, in our statement of the algorithm, we state that a step size for the multiplier vector is computed as
\bequation\label{eq:dual_update_beta}
  \beta_{y,k,j,l} \gets \arg\min_{\beta \in [0,1]} \|g_{k,j,l} + J_{k,j,l}^T(y_{k,j,l} + \beta \delta_{y,k,j,l})\|_2.
\eequation
On the other hand, for the slack multipliers, the algorithm first computes
\bequation\label{eq:z_step}
  \delta_{z,k,j,l} = \mu_{k,j} S_{k,j,l}^{-1} e - z_{k,j,l} - S_{k,j,l}^{-1} Z_{k,j,l} \dtilde_{s,k,j,l},
\eequation
then, to compute the step size $\beta_{z,k,j,l}$ satisfying the fraction-to-the-boundary rule
\bequationNN
  z_{k,j,l} + \beta_{z,k,j,l} \delta_{z,k,j,l} \geq (1 - \eta_1) z_{k,j,l},
\eequationNN
computes the limiting value
\bequation\label{eq:dual_update_beta_z}
  \beta_{z,k,j,l} = \min \left\{ 1, \min_{\left\{i :  \delta_{z,k,j,l}^{(i)} < 0 \right\}} \frac{-\eta_1 z^{(i)}_{k,j,l}}{\delta_{z,k,j,l}^{(i)}} \right\}.
\eequation

Finally, we state the termination conditions for the loop over index $l$ so that either the algorithm triggers a decrease in the barrier parameter or terminates completely so that the subsequent iteration of Algorithm~\ref{alg:ra} is reached.  Since the method ensures that $s_{k,j,l} > 0$ and $z_{k,j,l} > 0$ for all $(k,j,l)$, there is no need for the conditions to check positivity of the slack variables or their corresponding multipliers.  Therefore, one possibility for termination is that the algorithm verifies that \eqref{eq:sample:termination} holds with
\begin{multline}\label{eq.inplace}
  \left\{ \baligned (x_{k,j},s_{k,j},&y_{\Ecal,k,j},y_{\Ical,k,j},z_{k,j}) \\ &\text{or} \\(x_{k,j,l},s_{k,j,l},&y_{\Ecal,k,j,l},y_{\Ical,k,j,l},z_{k,j,l}) \ealigned \right\} \\ \text{in place of}\ (x_{k,M_k},s_{k,M_k},y_{\Ecal,k,M_k},y_{\Ical,k,M_k},z_{k,M_k}).
\end{multline}
On the other hand, the loop over $l$ terminates if an infeasible stationary point has been reached (consistent with the termination condition in Algorithm~\ref{alg:ra}) or, for some user-defined constants $(\epsilon_{\mu,1},\epsilon_{\mu,2},\epsilon_{\mu,3},\epsilon_{\mu,4},\epsilon_{\mu,5}) \in \R{}_{>0} \times \R{}_{>0} \times \R{}_{>0} \times \R{}_{>0} \times \R{}_{>0}$,
\begin{subequations}\label{eq:barrier_termination}
  \begin{align}
    \|\nabla f_{\Scal_k}(x_{k,j,l}) + \nabla c_\Ecal(x_{k,j,l})y_{\Ecal,k,j,l} + \nabla c_\Ical(x_{k,j,l})y_{\Ical,k,j,l}\|_{\infty} &\leq \epsilon_{\mu,1} \mu_{k,j}, \label{eq:barrier_termination_1} \\
    \|y_{\Ical,k,j,l} - z_{k,j,l}\|_\infty &\leq \epsilon_{\mu,2} \mu_{k,j} \label{eq:barrier_termination_2} \\
    \|c_\Ecal(x_{k,j,l})\|_{\infty} &\leq \epsilon_{\mu,3} \mu_{k,j}, \label{eq:barrier_termination_3} \\
    \|c_\Ical(x_{k,j,l}) + s_{k,j,l}\|_{\infty} &\leq \epsilon_{\mu,4} \mu_{k,j}, \label{eq:barrier_termination_4} \\
    \|S_{k,j,l} z_{k,j,l} - \mu_{k,j} e\|_{\infty} &\leq \epsilon_{\mu,5} \mu_{k,j}, \label{eq:barrier_termination_5}
  \end{align}
\end{subequations}
in which case a decrease in the barrier parameter is triggered.

\begin{algorithm}[H]
  \caption{Interior-Point Method}
  \label{alg.ipm}
  \footnotesize
  \begin{algorithmic}[1]
    \State \textbf{Inputs:} initial point $(x_{k,0},s_{k,0},y_{\Ecal,k,0},y_{\Ical,k,0},z_{k,0})$ and barrier parameter $\mu_{k,0} \in \R{}_{>0}$
    \State \textbf{Require:} $\omega \in \R{}_{>0}$, $\epsilon_v \in (0,1]$, $\kappa \in (0,1)$, $\psi \in \R{}_{>0}$, $\theta \in \R{}_{>0}$, $\zeta \in \R{}_{>0}$, $\sigma \in (0,1)$, $\epsilon_2 \in \R{}_{>0}$, $\epsilon_3 \in (0,1)$, $\tau \in (0,1)$, $\eta_1 \in (0,1)$, $\eta_2 \in (0,1)$, $\epsilon_\pi\in\R{}_{>0}$, $\{\epsilon_{\mu,a} \in \R{}_{>0}\}_{a\in[5]}$
    \For{$j = 0, 1, 2, \dots$}
      \If{$(x_{k,j}, s_{k,j}, y_{\Ecal,k,j}, y_{\Ical,k,j}, z_{k,j})$ satisfies \eqref{eq:sample:termination} with \eqref{eq.inplace}} \label{step:sub:termination}
        \State Set $M_k \gets j$ and \textbf{break}
      \EndIf
      \If{$(x_{k,j},s_{k,j})$ satisfies \eqref{eq:infeasible_stationary} and either $c_\Ecal(x_{k,j}) \neq 0$ or $c_\Ical(x_{k,j}) \not\leq 0$}
        \State Set $M_k \gets j$ and \textbf{break}
      \EndIf
      \State Set $(x_{k,j,0}, s_{k,j,0}, y_{\Ecal,k,j,0}, y_{\Ical,k,j,0}, z_{k,j,0}) \gets (x_{k,j}, s_{k,j}, y_{\Ecal,k,j}, y_{\Ical,k,j}, z_{k,j})$
      \State Choose $\pi_{k,j,-1} \in \R{}_{>0}$
      \For{$l = 0, 1, 2, \dots$}
        \If{$(x_{k,j,l}, s_{k,j,l}, y_{\Ecal,k,j,l}, y_{\Ical,k,j,l}, z_{k,j,l})$ satisfies \eqref{eq:barrier_termination} or \eqref{eq:sample:termination} with \eqref{eq.inplace}}
          \State Set $L_{k,j} = l$ and \textbf{break}
        \EndIf
        \If{$(x_{k,j,l},s_{k,j,l})$ satisfies \eqref{eq:infeasible_stationary} and either $c_\Ecal(x_{k,j,l}) \neq 0$ or $c_\Ical(x_{k,j,l}) \not\leq 0$}
          \State Set $L_{k,j} = l$ and \textbf{break}
        \EndIf
        \State Compute $v_{k,j,l}$ satisfying the Normal Step Condition
        \State Compute $(d_{k,j,l},\delta_{y,k,j,l})$ as an approximate solution of \eqref{eq.linsys}
        \While{$(d_{k,j,l},\delta_{y,k,j,l})$ does not satisfy any of Termination Tests 1, 2, or 3}
          \State Apply the Hessian Modification Strategy
          \State Compute an improved solution $(d_{k,j,l},\delta_{y,k,j,l})$ of \eqref{eq.linsys}
        \EndWhile
        \State Set $\pi_{k,j,l} \gets \pi_{k,j,l-1}$
        \If{$(d_{k,j,l},\delta_{y,k,j,l})$ does not satisfy TT1 but satisfies TT2}
          \State Set $d_{k,j,l} \gets 0$
        \ElsIf{$(d_{k,j,l},\delta_{y,k,j,l})$ satisfies neither TT1 nor TT2 (so it satisfies TT3)}
          \State Compute $\widehat{\pi}_{k,j,l}$ by \eqref{eq.pi_trial}
          \If{$\pi_{k,j,l} < \widehat{\pi}_{k,j,l}$}
            \State Set $\pi_{k,j,l} \gets \widehat{\pi}_{k,j,l} + \epsilon_\pi$
          \EndIf
        \EndIf
        \State Set $\dtilde_{k,j,l} \gets (d_{x,k,j,l}, S_{k,j,l} d_{s,k,j,l})$
        \State Set $\delta_{z,k,j,l}$ by \eqref{eq:z_step}
        \State Compute $\alpha_{k,j,l}^{\max}$ by \eqref{eq:alpha_max}
        \For{$i = 0, 1, 2,\dots$ until \eqref{eq:armijo} holds}
          \State Set $\alpha_{k,j,l} \gets 2^{-i}\alpha_{k,j,l}^{\max}$
        \EndFor
        \State Compute $\beta_{y,k,j,l}$ by \eqref{eq:dual_update_beta}
        \State Set $\beta_{z,k,j,l}$ by \eqref{eq:dual_update_beta_z}
        \State Set $x_{k,j,l+1} \gets x_{k,j,l} + \alpha_{k,j,l} \dtilde_{x,k,j,l}$ and $s_{k,j,l+1} \gets s_{k,j,l} + \alpha_{k,j,l} \dtilde_{s,k,j,l}$
        \State Set $s_{k,j,l+1} \gets \max\{s_{k,j,l+1}, - c_\Ical(x_{k,j,l+1})\}$
        \State Set $y_{\Ecal,k,j,l+1} \gets y_{\Ecal,k,j,l} + \beta_{y,k,j,l} \delta_{y_\Ecal,k,j,l}$ and $y_{\Ical,k,j,l+1} \gets y_{\Ical,k,j,l} + \beta_{y,k,j,l} \delta_{y_\Ical,k,j,l}$
        \State Set $z_{k,j,l+1} \gets z_{k,j,l} + \beta_{z,k,j,l} \delta_{z,k,j,l}$
      \EndFor
      \State Set
        $(x_{k,j+1}, s_{k,j+1}, y_{\Ecal,k,j+1}, y_{\Ical,k,j+1}, z_{k,j+1})
        \gets (x_{k,j,L_{k,j}}, s_{k,j,L_{k,j}}, y_{\Ecal,k,j,L_{k,j}}, y_{\Ical,k,j,L_{k,j}}, z_{k,j,L_{k,j}})$
      \If{the loop over $l$ terminated because \eqref{eq:sample:termination} holds with \eqref{eq.inplace}}
        \State Set $\mu_{k,j+1} \gets \mu_{k,j}$
      \Else
        \State Set $\mu_{k,j+1} \in (0,\mu_{k,j})$
      \EndIf
    \EndFor
  \end{algorithmic}
\end{algorithm}

\section{Convergence Analysis}\label{sec.analysis}

We begin our analysis by summarizing relevant results from \cite{CurtScheWaec2010}, applied to our setting of Algorithm~\ref{alg.ipm}, that give convergence guarantees for solving \eqref{eq:sampled_slack} for all generated $k \in \N{}$.  The following theorem covers potential finite termination of the algorithm, as well as potential convergence to infeasible stationary points or to points at which a constraint qualification fails to hold.

\btheorem\label{th:sub}
  Consider any generated $k \in \N{}$ for which Algorithm~\ref{alg:ra} calls Algorithm~\ref{alg.ipm}. Suppose that Assumptions~\ref{ass:full} and \ref{ass:sampled} hold, and that for any generated $j \in \N{}$ one has that (i) after a finite number of applications of the Hessian Modification Strategy, the matrix $W_{k,j,l}$ satisfies $W_{k,j,l} \succeq 2\theta I$ and is uniformly bounded in norm (by a constant that may depend on $k$, but not on $j$), and (ii) the iterative linear system solver can solve \eqref{eq.linsys} to arbitrary accuracy.  Moreover, suppose that if $j \to \infty$, then the barrier parameter sequence vanishes, i.e., $\{\mu_{k,j}\} \searrow 0$.  Then, one of the following occurs with respect to Algorithm~\ref{alg.ipm}.
  \benumerate
    \item[(a)] For some $(j,l)$, the primal point $(x_{k,j,l},s_{k,j,l})$ is an infeasible stationary point and the algorithm terminates finitely.
    \item[(b)] For some $j$ and all $l$, neither the set of conditions \eqref{eq:barrier_termination} nor \eqref{eq:sample:termination} with \eqref{eq.inplace} hold, and the algorithm generates iterates satisfying either
    \bequationNN
      \lim_{l \to \infty} \|J_{k,j,l}^Tc_{k,j,l}\| = 0\ \ \text{while}\ \ \|c_{k,j,l}\|_\infty \geq \min\{\epsilon_{\mu,3},\epsilon_{\mu,4}\} \mu_{k,j}\ \ \text{for all}\ \ l
    \eequationNN
    i.e., the iterates remain infeasible and \eqref{eq:infeasible_stationary} holds in a limiting sense, or
    \bequationNN
      \lim_{l \to \infty} \|J_{k,j,l}^Tc_{k,j,l}\| = 0\ \ \text{and}\ \ \lim_{l \to \infty} \pi_{k,j,l} = \infty,
    \eequationNN
    i.e., \eqref{eq:infeasible_stationary} holds in a limiting sense and the penalty parameter diverges.  This latter event cannot occur if the sequence of constraint Jacobians $\{J_{k,j,l}\}_{l=0}^\infty$ is bounded and all limit points have full row rank.
    \item[(c)] For some $j$, the termination conditions in line~\ref{step:sub:termination} hold, $s_{k,j} > 0$, and $z_{k,j} > 0$.
  \eenumerate
\etheorem
\bproof
  By construction of Algorithm~\ref{alg.ipm}, the conclusions follow as in \cite{CurtScheWaec2010}; in particular, see Assumptions~3.1 and 3.3 along with Theorems~3.4 and 3.13.
\eproof

Theorem~\ref{th:sub} shows that the solve of \eqref{eq:sampled_slack} for each generated $k \in \N{}$ either: terminates finitely with an infeasible stationary point; iterates indefinitely due to infeasible stationarity or lack of a constraint qualification; or terminates finitely with a point satisfying \eqref{eq:sample:termination}.  With respect to Algorithm~\ref{alg:ra}, the only interesting case to consider is when the latter event occurs for all generated $k \in \N{}$.  Moreover, if any iterate generated by Algorithm~\ref{alg:ra} is an infeasible stationary point, then there is also nothing left to prove.  Therefore, in essence, the only situation that remains interesting to consider is when Algorithm~\ref{alg:ra} generates an infinite sequence of iterates.

To proceed in our analysis of Algorithm~\ref{alg:ra}, we need to introduce some additional notation due to the fact that the randomness in the selection of the sample sets causes the resulting iterate sequence to be a stochastic process.  (We assume for the sake of simplicity that neither the linear system solver nor any other aspect of Algorithm~\ref{alg.ipm} involves stochastic quantities.  Our analysis could allow such possibilities with similar arguments, but for our purposes it is sufficient to presume that the stochastic nature of the algorithm is due solely to the sample sets.)  Let the possible sequences of iterations performed by Algorithm~\ref{alg:ra} define a probability space $(\Omega, \Gcal, \P)$.  Let $\Kcal$ be the event that Algorithm~\ref{alg:ra} generates an infinite sequence of iterates. Let $\Gcal_0$ be a $\sigma$-algebra defined by the initial conditions of Algorithm~\ref{alg:ra} and for all $k \in \N{}$ with $k > 0$ let $\Gcal_k$ denote the $\sigma$-algebra defined by the initial conditions and $\{\Scal_0,\dots,\Scal_{k-1}\}$.  We can now introduce the trace $\sigma$-algebra with $\Fcal_k := \Gcal_k \cap \Kcal$ for all $k \in \N{}$ defining the filtration $\{\Fcal_k\}$, which has $\Fcal_0 \subseteq \Fcal_1 \subseteq \cdots \subseteq \Gcal \cap \Kcal =: \Fcal$.

Let us emphasize that under this perspective, the sample sets, iterate sequence, steps computed by the algorithm, etc.~are random variables. \textit{For example, $x_k$ for each $k \in \N{}$ is a random variable that is $\Fcal_k$-measurable.} Employing this analytical framework, we make the following assumption about the gradients and their estimates.  The first inequality in the assumption is common in stochastic optimization methods; see, e.g., \cite{BottCurtNoce2018}.  The second inequality is similar to, e.g., \cite[eq.~(2.7)]{berahas2025retrospective}.

\bassumption\label{ass:indef}
  Let $(\Omega, \Gcal, \P)$, $\Kcal$, $\{\Gcal_k\}$, $\{\Fcal_k\}$, and $\Fcal$ be defined as above.  There exist positive constants $\chi_{\nabla f}$ and $\chi$ such that, for all $k \in \N{}$, one has
  \begin{align*}
    \E[\|\nabla f(x_k) - \nabla f_{\Scal_k}(x_k)\|_\infty^2 | \Fcal_k] &\leq \frac{\chi_{\nabla f} \|\nabla f(x_k)\|_\infty^2 + \chi}{|\Scal_k|} \\ \text{and}\ \ 
    \E[\|\nabla f(x_{k+1}) - \nabla f_{\Scal_k}(x_{k+1})\|_\infty^2 | \Fcal_k] &\leq \E\left[ \left. \frac{\chi_{\nabla f} \|\nabla f(x_{k+1})\|_\infty^2 + \chi}{|\Scal_k|} \right| \Fcal_k \right].
  \end{align*}
\eassumption

For the sake of brevity, let us define for all $k \in \N{}$ the primal-dual iterate
\bequationNN
  b_k := (x_k,s_k,y_{\Ecal,k},y_{\Ical,k},z_k).
\eequationNN
Let the quantities $b$, $b_{k,0}$, and $b_{k,M_k}$ be defined similarly following the notation of Algorithm~\ref{alg:ra}.  Further, let us define the residual functions $\{r_a\}_{a=1}^5$ and $r_{1,\Scal}$ by
\bequation\label{eq.residuals_all}
  \begin{aligned}
    r_1(b) &= \nabla f(x) + \nabla c_\Ecal(x) y_\Ecal + \nabla c_\Ical(x) y_\Ical, \\
    r_{1,\Scal}(b) &= \nabla f_\Scal(x) + \nabla c_\Ecal(x) y_\Ecal + \nabla c_\Ical(x) y_\Ical, \\
    r_2(b) &= y_\Ical - z, \\
    r_3(b) &= c_\Ecal(x), \\
    r_4(b) &= c_\Ical(x) + s, \\ \text{and}\ \ 
    r_5(b) &= Sz.
  \end{aligned}
\eequation

The following lemma bounds the residual norms during each iteration.

\blemma\label{lem:kkt_recursion}
  For all $k \in \N{}$, one has that
  \begin{align}
    \|r_1(b_{k,M_k})\|_\infty \leq&\ \gamma_{k,1} \|r_1(b_{k,0})\|_\infty + \epsilon_{k,1} \label{lem:recur_hard} \\
    &\ + \|\nabla f(x_{k,M_k}) - \nabla f_{\Scal_k}(x_{k,M_k})\|_\infty \nonumber \\
    &\ + \gamma_{k,1} \|\nabla f(x_{k,0}) - \nabla f_{\Scal_k}(x_{k,0})\|_\infty, \nonumber
  \end{align}
  whereas for all $a \in \{2,3,4,5\}$ one has that
  \bequation\label{lem:recur_easy}
    \|r_a(b_{k,M_k})\|_\infty \leq \gamma_{k,a} \|r_a(b_{k,0})\|_\infty + \epsilon_{k,a}.
  \eequation
\elemma
\bproof
  Consider arbitrary $k \in \N{}$.  The inequalities in \eqref{lem:recur_easy} follow directly from \eqref{eq:sample:termination_2}--\eqref{eq:sample:termination_5}.  As for \eqref{lem:recur_hard}, one finds by the triangle inequality and \eqref{eq:sample:termination_1} that
  \begin{align*}
    \|r_1(b_{k,M_k})\|_\infty
    \leq&\ \|r_{1,\Scal_k}(b_{k,M_k})\|_\infty + \|r_1(b_{k,M_k}) - r_{1,\Scal_k}(b_{k,M_k})\|_\infty \\
    \leq&\ \gamma_{k,1} \|r_{1,\Scal_k}(b_{k,0})\|_\infty + \epsilon_{k,1} + \|r_1(b_{k,M_k}) - r_{1,\Scal_k}(b_{k,M_k})\|_\infty \\
    \leq&\ \gamma_{k,1} (\|r_1(b_{k,0})\|_\infty + \|r_{1,\Scal_k}(b_{k,0}) - r_1(b_{k,0})\|_\infty) + \epsilon_{k,1} \\
    &\ + \|r_1(b_{k,M_k}) - r_{1,\Scal_k}(b_{k,M_k})\|_\infty.
  \end{align*}
  Since $r_1(b) - r_{1,\Scal}(b) = \nabla f(x) - \nabla f_\Scal(x)$ for any $b$ and $\Scal$, the proof is complete.
\eproof

Our next lemma considers the effect of the perturbation of the slack variables and their multipliers that can occur in line~\ref{step.slack_perturbation} of Algorithm~\ref{alg:ra}.  (Note that the perturbation of the barrier parameter in this line of the algorithm has no effect on these quantities.)

\blemma\label{lem:perturbation}
  For all $k \in \N{}$, one has that
  \begin{align*}
    \|r_2(b_{k+1})\|_\infty &\leq \|r_2(b_{k,M_k})\|_\infty + \nu_{z,k}, \\
    \|r_4(b_{k+1})\|_\infty &\leq \|r_4(b_{k,M_k})\|_\infty + \nu_{s,k}, \\ \text{and}\ \ 
    \|r_5(b_{k+1})\|_\infty &\leq \|r_5(b_{k,M_k})\|_\infty + \|s_{k,M_k}\|_\infty \nu_{z,k} + \nu_{s,k} \|z_{k,M_k}\|_\infty + \nu_{s,k} \nu_{z,k},
  \end{align*}
  whereas for each $a \in \{1,3\}$ one has that $r_a(b_{k+1}) = r_a(b_{k,M_k})$.
\elemma
\bproof
  The only differences between $b_{k+1}$ and $b_{k,M_k}$ occur in the slack variables and their corresponding dual variables, which have no effect on the $r_1$ and $r_3$ residual values.  Hence, the equations for $a \in \{1,3\}$ hold trivially.  The bounds for the residuals $r_2$ and $r_4$ follow by the triangle inequality, $s_{k+1} \gets s_{k,M_k} + \iota_{s,k}$, $z_{k+1} \gets z_{k,M_k} + \iota_{z,k}$, and the fact that $\iota_{s,k}$ and $\iota_{z,k}$ are chosen to satisfy $\|\iota_{s,k}\|_\infty \leq \nu_{s,k}$ and $\|\iota_{z,k}\|_\infty \leq \nu_{z,k}$, respectively.  The bound for the residual $r_5$ follows from the triangle inequality, the fact that $r_5(b_{k+1}) - r_5(b_{k,M_k}) = S_{k,M_k} \iota_{z,k} + Z_{k,M_k} \iota_{s,k} + \iota_{s,k} \cdot \iota_{z,k}$, the facts that $\|\iota_{s,k}\|_\infty \leq \nu_{s,k}$ and $\|\iota_{z,k}\|_\infty \leq \nu_{z,k}$, and the fact that for two vectors of the same dimension with nonnegative entries, the $\ell_\infty$-norm of their componentwise-product is less than or equal to the product of the $\ell_\infty$-norms of the two vectors.
\eproof

Our main result employs the following technical lemma.

\blemma\label{lem:toeplitz}
  If $q_k \in \R{}_{\geq 0}$ for all $k \in \N{}$ has $\{q_k\} \to 0$ and $\gamma \in (0,1)$, then
  \bequationNN
    \left\{\sum_{k=0}^K q_k \gamma^{K-k}\right\} \to 0\ \ \text{as}\ \ K \to \infty.
  \eequationNN
\elemma
\bproof
  Since $\{q_k\} \to 0$, it follows that $\qbar := \max \{q_k : k \in \N{}\} \in\R{}_{\geq 0}$ exists and is finite. Now consider arbitrary $\underline{q} \in \R{}_{>0}$. Since $\{q_k\} \to 0$, there exists $K_{\underline{q}} \in \N{}$ such that $q_k \leq \underline{q}$ for all $k \geq K_{\underline{q}}$.  For any $K \in \N{}$ with $K > K_{\underline{q}}$, one finds that
  \begin{align}
    \sum_{k=0}^K q_k \gamma^{K-k}
    &= \sum_{k=0}^{K_{\underline{q}}} q_k \gamma^{K-k} + \sum_{k=K_{\underline{q}}+1}^K q_k \gamma^{K-k} \nonumber \\
    &\leq \qbar \sum_{k=0}^{K_{\underline{q}}} \gamma^{K-k} + \underline{q} \sum_{k=K_{\underline{q}}+1}^K \gamma^{K-k} \nonumber \\
    &\leq \qbar \frac{\gamma^{K-K_{\underline{q}}}}{1 - \gamma} + \frac{\underline{q}}{1 - \gamma}. \label{eq:toeplitz}
  \end{align}
  For any $\varepsilon \in \R{}_{>0}$, one can choose $\underline{q}$ sufficiently small relative to $1-\gamma$, then choose $K_{\underline{q}}$ such that $q_k \leq \underline{q}$ for all $k \geq K_{\underline{q}}$, and finally choose $K$ sufficiently large relative to $K_{\underline{q}}$ such that \eqref{eq:toeplitz} is less than $\varepsilon$.  Thus, the desired conclusion holds.
\eproof

We are now prepared to prove a convergence guarantee for Algorithm~\ref{alg:ra}.

\btheorem\label{th:main}
  Suppose that Assumptions~\ref{ass:full}, \ref{ass:sampled}, and \ref{ass:indef} hold and $\{|\Scal_k|\} \to \infty$.  Suppose that there exists $\gamma_a \in (0,1)$ such that $\gamma_{k,a} \leq \gamma_a$ for all $(k,a) \in \N{} \times [5]$, and that $\{\epsilon_{k,a}\}_{k \in \N{}} \to 0$ for all $a \in [5]$, $\{\nu_{s,k}\}_{k\in\N{}} \to 0$, and $\{\nu_{z,k}\}_{k\in\N{}} \to 0$.  Finally, suppose that there exist $(\kappa_s, \kappa_z) \in \R{}_{>0} \times \R{}_{>0}$ and $k_0 \in \N{}$ such that $\|s_{k,M_k}\|_\infty \leq \kappa_s$ and $\|z_{k,M_k}\|_\infty \leq \kappa_z$ for all $k \geq k_0$.  Then, the sequences of expected residuals satisfy
  \bequationNN
    \lim_{k \to \infty} \E[\|r_a(b_k)\|_\infty] = 0\ \ \text{for all}\ \ a \in [5].
  \eequationNN
\etheorem
\bproof
  Throughout the proof, let us define the constant
  \bequationNN
    \kappa_\chi := \sqrt{\chi_{\nabla f} \kappa_{\nabla f}^2 + \chi} \in \R{}_{>0}.
  \eequationNN

  Let us begin by considering the residual $r_1$.  Consider arbitrary $k \in \N{}$.  By Lemma~\ref{lem:perturbation}, one has $r_1(b_{k+1}) = r_1(b_{k,M_k})$, so Lemma~\ref{lem:kkt_recursion} along with the facts that $b_{k,0} = b_k$ and $x_{k,M_k} = x_{k+1}$ yields
  \begin{align*}
    \|r_1(b_{k+1})\|_\infty
    \leq&\ \gamma_{k,1} \|r_1(b_k)\|_\infty + \epsilon_{k,1} \\
    &\ + \|\nabla f(x_{k+1}) - \nabla f_{\Scal_k}(x_{k+1})\|_\infty + \gamma_{k,1} \|\nabla f(x_k) - \nabla f_{\Scal_k}(x_k)\|_\infty.
  \end{align*}
  Since $b_k$ is $\Fcal_k$-measurable and all terms are nonnegative, one obtains that
  \begin{align}
    \E[\|r_1(b_{k+1})\|_\infty | \Fcal_k]
    \leq&\ \gamma_{k,1} \|r_1(b_k)\|_\infty + \epsilon_{k,1} \nonumber \\
    &\ + \E[\|\nabla f(x_{k+1}) - \nabla f_{\Scal_k}(x_{k+1})\|_\infty | \Fcal_k] \nonumber \\
    &\ + \gamma_{k,1} \E[\|\nabla f(x_k) - \nabla f_{\Scal_k}(x_k)\|_\infty | \Fcal_k]. \label{eq:main_cond}
  \end{align}
  Let us now bound the last two terms on the right-hand side of \eqref{eq:main_cond}.  By Jensen's inequality (applied with concave $t \mapsto \sqrt{t}$), the first inequality in Assumption~\ref{ass:indef}, and the bound $\|\nabla f(x_k)\|_\infty \leq \|\nabla f(x_k)\|_2 \leq \kappa_{\nabla f}$ from Assumption~\ref{ass:full}, one finds
  \begin{align*}
    \E[\|\nabla f(x_k) - \nabla f_{\Scal_k}(x_k)\|_\infty | \Fcal_k]
    &\leq \( \E[\|\nabla f(x_k) - \nabla f_{\Scal_k}(x_k)\|_\infty^2 | \Fcal_k] \)^{1/2} \\
    &\leq \( \frac{\chi_{\nabla f} \|\nabla f(x_k)\|_\infty^2 + \chi}{|\Scal_k|} \)^{1/2}
     \leq \frac{\kappa_\chi}{\sqrt{|\Scal_k|}}.
  \end{align*}
  Similarly, by the second inequality in Assumption~\ref{ass:indef} and since $\|\nabla f(x_{k+1})\|_\infty \leq \kappa_{\nabla f}$ holds surely by Assumption~\ref{ass:full}, one finds that
  \begin{align*}
    &\ \E[\|\nabla f(x_{k+1}) - \nabla f_{\Scal_k}(x_{k+1})\|_\infty | \Fcal_k] \\
    \leq&\ \( \E[\|\nabla f(x_{k+1}) - \nabla f_{\Scal_k}(x_{k+1})\|_\infty^2 | \Fcal_k] \)^{1/2} \\
    \leq&\ \( \E\left[ \left. \frac{\chi_{\nabla f} \|\nabla f(x_{k+1})\|_\infty^2 + \chi}{|\Scal_k|} \right| \Fcal_k \right] \)^{1/2}
     \leq \frac{\kappa_\chi}{\sqrt{|\Scal_k|}}.
  \end{align*}
  Combining these bounds with \eqref{eq:main_cond} and taking total expectation, it follows that
  \bequation\label{eq:main_recur}
    \E[\|r_1(b_{k+1})\|_\infty] \leq \gamma_{k,1} \E[\|r_1(b_k)\|_\infty] + q_k,\ \ \text{where}\ \ q_k := \epsilon_{k,1} + (1 + \gamma_{k,1}) \frac{\kappa_\chi}{\sqrt{|\Scal_k|}}.
  \eequation
  Since $k \in \N{}$ was arbitrary, \eqref{eq:main_recur} holds for all $k \in \N{}$.  Applying it recursively, one finds for any $K \in \N{}$, with the convention that $\prod_{i=K+1}^K \gamma_{i,1} := 1$, that
  \begin{align}
    \E[\|r_1(b_{K+1})\|_\infty]
    &\leq \( \prod_{k=0}^K \gamma_{k,1} \) \|r_1(b_0)\|_\infty + \sum_{k=0}^K \( \prod_{i=k+1}^K \gamma_{i,1} \) q_k \nonumber \\
    &\leq \gamma_1^{K+1} \|r_1(b_0)\|_\infty + \sum_{k=0}^K q_k \gamma_1^{K-k}, \label{eq:main_unrolled}
  \end{align}
  where the second inequality follows since $\gamma_{k,1} \leq \gamma_1$ and $q_k \geq 0$ for all $k \in \N{}$.  The first term on the right-hand side of \eqref{eq:main_unrolled} vanishes as $K \to \infty$ since $\gamma_1 \in (0,1)$ and $b_0$ is fixed.  As for the second term, since $\{\epsilon_{k,1}\} \to 0$, $\{|\Scal_k|\} \to \infty$, and $\gamma_{k,1} \leq 1$ for all $k \in \N{}$, one has $\{q_k\} \to 0$, so this term vanishes as $K \to \infty$ by Lemma~\ref{lem:toeplitz}.  This proves the desired conclusion for $a = 1$.

  Now consider $a \in \{2,3,4,5\}$.  Combining Lemmas~\ref{lem:kkt_recursion} and \ref{lem:perturbation} with the facts that $b_{k,0} = b_k$ for all $k \in \N{}$ and that $\|s_{k,M_k}\|_\infty \leq \kappa_s$ and $\|z_{k,M_k}\|_\infty \leq \kappa_z$ for all $k \geq k_0$,
  \begin{align*}
    \|r_2(b_{k+1})\|_\infty \leq&\ \gamma_{k,2} \|r_2(b_k)\|_\infty + \epsilon_{k,2} + \nu_{z,k}\ \text{for all}\ k \in \N{}, \\
    \|r_3(b_{k+1})\|_\infty \leq&\ \gamma_{k,3} \|r_3(b_k)\|_\infty + \epsilon_{k,3}\ \text{for all}\ k \in \N{}, \\
    \|r_4(b_{k+1})\|_\infty \leq&\ \gamma_{k,4} \|r_4(b_k)\|_\infty + \epsilon_{k,4} + \nu_{s,k}\ \text{for all}\ k \in \N{}, \\ \text{and}\ \ 
    \|r_5(b_{k+1})\|_\infty \leq&\ \gamma_{k,5} \|r_5(b_k)\|_\infty + \epsilon_{k,5} + \kappa_s \nu_{z,k} + \kappa_z \nu_{s,k} + \nu_{s,k} \nu_{z,k}\ \text{for all}\ k \geq k_0.
  \end{align*}
  Unlike \eqref{eq:main_recur}, these inequalities hold surely, i.e., for every realization of the sample sets.  Each has the form
  \bequation\label{eq:main_recur_a}
    \|r_a(b_{k+1})\|_\infty \leq \gamma_{k,a} \|r_a(b_k)\|_\infty + q_{a,k}\ \ \text{for all}\ \ k \geq k_a,
  \eequation
  where $\{q_{a,k}\}$ is a deterministic nonnegative sequence with $\{q_{a,k}\} \to 0$ under the conditions of the theorem, and where $k_a := 0$ for $a \in \{2,3,4\}$ and $k_5 := k_0 + 1$.  Moreover, the residual at index $k_a$ is bounded by a deterministic constant $\kappa_{r_a} \in \R{}_{\geq 0}$: for $a \in \{2,3,4\}$ one can take $\kappa_{r_a} := \|r_a(b_0)\|_\infty$ since $b_0$ is fixed, whereas for $a = 5$ it follows from line~\ref{step.slack_perturbation} of Algorithm~\ref{alg:ra} that $\|s_{k_0+1}\|_\infty \leq \kappa_s + \nu_{s,k_0}$ and $\|z_{k_0+1}\|_\infty \leq \kappa_z + \nu_{z,k_0}$, so that one can take $\kappa_{r_5} := (\kappa_s + \nu_{s,k_0})(\kappa_z + \nu_{z,k_0}) \geq \|S_{k_0+1} z_{k_0+1}\|_\infty = \|r_5(b_{k_0+1})\|_\infty$.  Hence, applying \eqref{eq:main_recur_a} recursively as in \eqref{eq:main_unrolled}, one finds surely for any $K \in \N{}$ with $K \geq k_a$ that
  \begin{align*}
    \|r_a(b_{K+1})\|_\infty
    &\leq \gamma_a^{K+1-k_a} \kappa_{r_a} + \sum_{k=k_a}^K q_{a,k} \gamma_a^{K-k} \\
    &= \gamma_a^{K+1-k_a} \kappa_{r_a} + \sum_{m=0}^{K-k_a} q_{a,m+k_a} \gamma_a^{(K-k_a)-m}.
  \end{align*}
  The right-hand side is deterministic, so the same bound holds for $\E[\|r_a(b_{K+1})\|_\infty]$.  The first term vanishes as $K \to \infty$ since $\gamma_a \in (0,1)$, and the second term vanishes by Lemma~\ref{lem:toeplitz} applied to the nonnegative sequence $\{q_{a,m+k_a}\}_{m \in \N{}}$, which converges to zero.  This proves the desired conclusion for $a \in \{2,3,4,5\}$.
\eproof

\section{Implementation}\label{sec.implementation}

The purpose of this section is to describe the details of an implementation of the algorithm proposed in \S\ref{sec.algorithm}, namely, that presented as Algorithms~\ref{alg:ra} and~\ref{alg.ipm}. Let us refer to our particular implementation as the retrospective approximation interior-point method (RA-IPM). The details to be described are the techniques for prescaling, choosing the sample sizes $\{N_k\}$, selecting the perturbations $\{(\iota_{s,k},\iota_{z,k},\iota_{\mu,k})\}$, computing the normal steps $\{v_{k,j,l}\}$, computing the approximate solutions of \eqref{eq.linsys}, and performing the line searches and iterate updates. The values of the input parameters employed are provided in Tables~\ref{tab:parameters} and~\ref{tab:parameters2} at the end of this section. The software environment used for our experiments is described in \S\ref{sec.numerical}.

\subsection{Prescaling}\label{sec.prescaling}

During initialization, a problem is prescaled by the gradient-based strategy employed in Ipopt~\cite{wachter2006implementation}. For example, the objective is multiplied by $\max\{s^{\rm scale}_{\min}, \min\{1, g_{\max}/\|\nabla f(x_0)\|_\infty\}\}$, where $g_{\max}$ and $s^{\rm scale}_{\min}$ are user-defined. Also, each constraint function is multiplied by an analogous factor determined by the corresponding row of the constraint Jacobians (i.e., $\nabla c_\Ecal$ and $\nabla c_\Ical$) at $x_0$.  These factors are held fixed throughout the solve. (Also, the same factors are supplied to every method compared in \S\ref{sec.numerical}.)  All tolerances reported below refer to the scaled problem.

\subsection{Sampling}\label{sec.impl_sampling}

Algorithm~\ref{alg:ra} requires a sequence of samples $\{\Scal_k\}$ with sample sizes $\{N_k\}$. Rather than prescribe the sequence of sample sizes in advance, our implementation generates it adaptively following the variance-based rules proposed in~\cite{berahas2025retrospective,bollapragada2018adaptive}. We describe this procedure in this subsection. Throughout, we presume that the optimization problem of interest has a full set $\Scal_{\rm full}$ of $N$ realizations of $\xi$, which is generated prior to the first iteration (in the same manner for all methods in~\S\ref{sec.numerical}.)  Every sample employed by the implementation is therefore a subset of $\Scal_{\rm full}$ that is drawn without replacement, subject to an acceptance rule explained in~\S\ref{sec.numerical}.

Let us now describe the procedure for selecting the sample sizes. For any sample~$\Scal$, let $r_\Scal(b)$ denote the stacked vector of residuals as defined in \eqref{eq.residuals_all} (using $r_{1,\Scal}$, not~$r_1$).  In iteration $k=0$, an initial sample size $N_0$ is employed.  Otherwise, for $k > 0$, at the beginning of the iteration, a subset $\widetilde\Scal_{k-1} \subseteq \Scal_{\rm full}$ of size $\widetilde N_{k-1} = N_{k-1}$ is drawn.  This draw is made independently of $\Scal_{k-1}$, since $x_{k,0}$ inherits information from the inner solve performed over $\Scal_{k-1}$, meaning that a variance estimated over $\Scal_{k-1}$ would be optimistic~\cite{berahas2025retrospective}.  Over $\widetilde\Scal_{k-1}$, the implementation computes the sample variance
\bequation\label{eq:Z_V_estimates}
  \operatorname{Var}_{\widetilde\Scal_{k-1}}(x_{k,0}) := \frac{1}{\widetilde N_{k-1}} \sum_{\xi \in \widetilde\Scal_{k-1}} \big\| \nabla F(x_{k,0},\xi) - \nabla f_{\widetilde\Scal_{k-1}}(x_{k,0}) \big\|_2^2,
\eequation
and forms the sample size for which the standard deviation of the sampled gradient would be a fraction $\tilde\theta \in (0,1)$ of $\|r_{\widetilde\Scal_{k-1}}(b_{k,0})\|_2$, namely, the target
\bequation\label{eq:N_k_var}
  N^{\rm var}_k := \left\lceil \frac{\operatorname{Var}_{\widetilde\Scal_{k-1}}(x_{k,0})} {\tilde\theta^2 \|r_{\widetilde\Scal_{k-1}}(b_{k,0})\|_2^2} \right\rceil.
\eequation
The target $N^{\rm var}_k$ is determined entirely by the quantities estimated
over $\widetilde\Scal_{k-1}$, and as such may fall below $N_{k-1}$.  A floor and
a ceiling are therefore imposed, which yields
\bequation\label{eq:N_k_rule}
  N_k := \min\big\{ N,\ \hat\beta N_{k-1},\
    \max\{ 2N_{k-1},\ N^{\rm var}_k \} \big\},
\eequation
where $\hat\beta \geq 2$ is a user-defined constant. Note that the presence of $2N_{k-1}$ in the max forces the sample size to at least double in each outer iteration so that the full sample is reached within a budget, call it $\bar K$, of iterations; see Tables~\ref{tab:parameters}--\ref{tab:parameters2}.  The ceiling of $\hat\beta N_{k-1}$, on the other hand, prevents a single large variance estimate from causing the sample to reach the full sample set $\Scal_{\rm full}$ in a single large jump in sample size.

The set $\Scal_k$ is then formed by augmenting $\widetilde\Scal_{k-1}$ with $N_k - \widetilde N_{k-1}$ further realizations drawn from $\Scal_{\rm full} \setminus \widetilde\Scal_{k-1}$ so that $\widetilde\Scal_{k-1} \subset \Scal_k$.  The gradient over $\Scal_k$ is recomputed rather than assembled from \eqref{eq:Z_V_estimates}.  Once $N_k = N$ is attained, the sampled problem coincides with \eqref{eq:full} and the outer loop has no further sample to choose. At this time, the relative tolerances $\{\gamma_{k,a}\}$ in \eqref{eq:sample:termination} are set to zero, and what remains is a single solve of~\eqref{eq:full} to the absolute tolerances.  We remark that, $\Scal_{\rm full}$ being finite, this occurs after finitely many outer iterations rather than in the limit presumed in \S\ref{sec.analysis}.

\subsection{Warm-starting}\label{sec.impl_warmstart}

Algorithm~\ref{alg:ra} permits the slack variables, their corresponding dual variables, and the barrier parameter to be increased before the subsequent sampled problem is solved.  Consecutive sampled problems differ only through their samples, so their solutions might lie close to one another, meaning that a warm start of Algorithm~\ref{alg.ipm} could be effective.  That said, a barrier parameter carried forward may nonetheless be too small for a central path that has moved, and slack variables and/or their duals close to zero by the previous barrier subproblem might leave the next one ill-conditioned at a point that is no longer nearly optimal.  Whether the previous solution is worth reusing thus depends on how far the sampled objective has moved.  Once $\Scal_k$ has been drawn, we compute the relative gradient drift
\bequation\label{eq:outer_restart_rel_drift}
  \rho_k := \frac{\|\nabla f_{\Scal_k}(x_{k,0})
                   - \nabla f_{\Scal_{k-1}}(x_{k,0})\|_\infty}
                 {\max\{\|\nabla f_{\Scal_{k-1}}(x_{k,0})\|_\infty,\
                         \underline\rho\}},
\eequation
where $\underline\rho \in \R{}_{>0}$ safeguards the denominator.  (With the value $\underline\rho = 1$ in Table~\ref{tab:parameters}, $\rho_k$ measures the drift relative to the larger of the gradient norm and one, so that the test is on the absolute drift when $\|\nabla f_{\Scal_{k-1}}(x_{k,0})\|_\infty \leq 1$.)  A small value of $\rho_k$ suggests that the central path of the new sampled problem lies close (at least locally) to the one that the inner loop has just followed, whereas a large value suggests that it is not.  For a user-defined threshold $\bar\rho \in \R{}_{>0}$, we say that the transition into outer iteration $k \geq 1$ is \emph{warm-started} when $\rho_k \leq \bar\rho$ and the preceding call to Algorithm~\ref{alg.ipm} terminated through \eqref{eq:sample:termination}. The perturbations $\iota_{s,k-1}$ and $\iota_{z,k-1}$ chosen at the end of outer iteration $k-1$ and the barrier parameter $\mu_{k,0}$ are both governed by this test.  The primal-dual iterate $(x_{k,M_k},s_{k,M_k},y_{\Ecal,k,M_k},y_{\Ical,k,M_k},z_{k,M_k})$ is carried forward in every case, as line~\ref{step.slack_perturbation} of Algorithm~\ref{alg:ra} prescribes. A transition that is not warm-started preserves only the barrier parameter and the barrier-stage budget.  Other strategies are possible. The one described here worked well in our tests.

Within Algorithm~\ref{alg.ipm}, our implementation decreases the barrier parameter by the monotone rule employed in Ipopt~\cite{wachter2006implementation}, namely,
\bequation\label{eq:mu_update}
  \mu_{k,j+1} \gets \max\{\underline\mu, \min\{\kappa_\mu \mu_{k,j}, \mu_{k,j}^{\theta_\mu}\}\}
\eequation
for user-defined $\kappa_\mu \in (0,1)$, $\theta_\mu \in (1,2)$, and $\underline\mu \in \R{}_{>0}$, and it performs at most $\bar J$ barrier stages, i.e., iterations of the loop over $j$, per sampled problem.

For the barrier parameter, for a user-defined $c_\mu \geq 1$, the implementation sets
\bequation\label{eq:mu_hotstart}
  \mu_{k,0} \gets
  \begin{cases}
    \min\{c_\mu \mu_{k-1,M_{k-1}}, \mu_0\},
      & \text{if outer iteration $k$ is warm-started}, \\[2pt]
    \mu_0, & \text{otherwise},
  \end{cases}
\eequation
where $\mu_{k-1,M_{k-1}}$ denotes the barrier parameter upon termination of outer iteration $k-1$.  Both cases are permitted since $\iota_{\mu,k-1} \in [0,\infty)$ is required there and the barrier parameter is nonincreasing within Algorithm~\ref{alg.ipm}, so $\mu_{k-1,M_{k-1}} \leq \mu_0$.  The factor~$c_\mu$ loosens the carried value, since the central path has moved even when the two samples are close. The cap at $\mu_0$ keeps the carried value no larger than the one with which outer iteration $k=0$ begins.  A warm-started iteration thus begins at a barrier parameter no larger than $\mu_0$, so the $\bar J$ updates available to Algorithm~\ref{alg.ipm} would carry it below the value that $\bar J$ updates from $\mu_0$ would reach.  The number of stages available in outer iteration~$k$ is therefore reduced by the number of applications of \eqref{eq:mu_update} that would be needed to bring $\mu_0$ to or below $\mu_{k,0}$, so that the two trajectories terminate at comparable values of the barrier parameter.

As for the slack variables, a perturbation of $s$ alone leaves $x$, $y_\Ecal$, $y_\Ical$, and $z$ unchanged, so by Lemma~\ref{lem:perturbation} its only effect is on the residuals $r_4$ and $r_5$, to which it contributes $\iota_{s,k}$ and $Z_{k,M_k}\iota_{s,k}$, respectively.  Requiring each contribution to consume at most a fraction $c_s \in \R{}_{>0}$ of the corresponding tolerance yields the budget
\bequation\label{eq:slack_bump_pre}
  \tilde\epsilon_{s,k} := c_s \min\left\{ \epsilon_{k,4},\
    \frac{\epsilon_{k,5}}{\|z_{k,M_k}\|_\infty} \right\}.
\eequation
Letting $\underline s_k := \min_{i \in [q]} s^{(i)}_{k,M_k}$ denote the smallest slack variable returned by the inner solve, we raise to the level $\underline s_k + \tilde\epsilon_{s,k}$ every slack variable lying below it, i.e., we set
\bequation\label{eq:slack_bump}
  \iota^{(i)}_{s,k} :=
  \begin{cases}
    \max\{ \underline s_k + \tilde\epsilon_{s,k} - s^{(i)}_{k,M_k},\ 0 \},
      & \text{if outer iteration $k$ is warm-started}, \\[4pt]
    0, & \text{otherwise},
  \end{cases}
\eequation
for all $i \in [q]$; equivalently, $s^{(i)}_{k+1}$ is the larger of $s^{(i)}_{k,M_k}$ and this floor.  The floor is thus set relative to the smallest slack variable rather than at an absolute threshold, so that on a warm-started iteration the component attaining $\underline s_k$ is raised by exactly $\tilde\epsilon_{s,k}$, no component is raised by more, and a component already above the floor is left unchanged.  One therefore has
\begin{align*}
  \|\iota_{s,k}\|_\infty &= \tilde\epsilon_{s,k} \leq c_s\epsilon_{k,4} \\
  \text{and}\ \ \|Z_{k,M_k}\iota_{s,k}\|_\infty &\leq \|z_{k,M_k}\|_\infty \|\iota_{s,k}\|_\infty \leq c_s\epsilon_{k,5},
\end{align*}
which are the conditions required in line~\ref{step.slack_perturbation} of Algorithm~\ref{alg:ra}.  The condition in the first line holds with equality, whereas the second line may be strict, since the largest component of $z_{k,M_k}$ and the largest component of $\iota_{s,k}$ need not be the same.  Finally, the multipliers corresponding to the bounds on the slack variables are left unchanged.

\subsection{Normal Step Computation}\label{sec.normal_impl_step}

The normal step is computed by the dogleg procedure stated as Algorithm~\ref{alg:dogleg_normal}, in which the Newton step is obtained by solving the normal equations $J_{k,j,l}J_{k,j,l}^T\varpi_{k,j,l} = -c_{k,j,l}$ using a matrix factorization.  Should the factorization fail, as it might for rank-deficient $J_{k,j,l}$, a positive-definite term $\delta I$ is added to $J_{k,j,l}J_{k,j,l}^T$ and increased by a factor $\varsigma_v$ until it succeeds. In such a case, the Newton step solves a regularized problem, and on its own need not reduce linearized infeasibility by as much as the Cauchy point. The step returned is therefore taken from the segment joining the two, over which the Cauchy point is attained at $a^\star = 0$, so that the Normal Step Condition of \S\ref{sec.interior} holds with $\epsilon_v = 1$.  The scaling $\omega$ is increased by a factor of ten, to at most $\omega_{\max}$, whenever a full step is accepted whose normal component lies on the trust-region boundary, and is reset at the start of each sampled problem \cite{CurtHubeScheWaec2012}.

\begin{algorithm}[ht]
  \caption{Normal Step Computation}
  \label{alg:dogleg_normal}
  \footnotesize
  \begin{algorithmic}[1]
    \State \textbf{Given} $\omega \in \R{}_{>0}$, $\delta_v \in \R{}_{>0}$ and $\varsigma_v \in \R{}_{>1}$
    \State \textbf{if} $J_{k,j,l}^Tc_{k,j,l} = 0$ \textbf{then return} $v_{k,j,l} \gets 0$
    \State Set $\Omega_{k,j,l} \gets \omega\|J_{k,j,l}^Tc_{k,j,l}\|_2$
    \State Set the Cauchy point
      $v^C_{k,j,l} \gets -\min\left\{\omega,\
      \dfrac{\|J_{k,j,l}^Tc_{k,j,l}\|_2^2}{\|J_{k,j,l}J_{k,j,l}^Tc_{k,j,l}\|_2^2}\right\}
      J_{k,j,l}^Tc_{k,j,l}$
    \State Set $\delta \gets 0$
    \Loop
      \State Attempt to factor $J_{k,j,l}J_{k,j,l}^T + \delta I$ and solve
        $(J_{k,j,l}J_{k,j,l}^T + \delta I)\varpi_{k,j,l} = -c_{k,j,l}$
      \State \textbf{if} successful \textbf{then break}; \textbf{else}
        set $\delta \gets \delta_v$ if $\delta = 0$ and $\delta \gets \varsigma_v\delta$ otherwise
    \EndLoop
    \State Set $v^N_{k,j,l} \gets J_{k,j,l}^T\varpi_{k,j,l}$
    \If{$\|v^C_{k,j,l}\|_2 = \Omega_{k,j,l}$}
      \State \textbf{return} $v_{k,j,l} \gets v^C_{k,j,l}$
    \Else
      \State Compute $a^\star$ as a solution of
      \bequationNN
        \min_{a \in [0,1]}\ \|c_{k,j,l} + J_{k,j,l}(v^C_{k,j,l} + a(v^N_{k,j,l} - v^C_{k,j,l}))\|_2\ \ \st\ \ \|v^C_{k,j,l} + a(v^N_{k,j,l} - v^C_{k,j,l})\|_2 \leq \Omega_{k,j,l}
      \eequationNN
      \State \textbf{return} $v_{k,j,l} \gets v^C_{k,j,l} + a^\star(v^N_{k,j,l} - v^C_{k,j,l})$
    \EndIf
  \end{algorithmic}
\end{algorithm}

\subsection{Primal-Dual Step Computation}\label{sec.impl_step}

Algorithm~\ref{alg.ipm} requires $(d_{k,j,l},\delta_{y,k,j,l})$ to be an approximate solution of \eqref{eq.linsys} satisfying one of Termination Tests~1, 2, or~3 of~\S\ref{sec.interior}.  Our implementation computes such a step with an iterative method, described below, that ultimately tests every iterate and accepts the first one that passes a termination test.  As preconditioner we employ the constraint preconditioner obtained from \eqref{eq.linsys} by replacing the leading block with $D_{k,j,l} := \diag(\max\{|\diag(W_{k,j,l})|, \delta_P\})$; see~\cite{KellGoulWath2000}.  Its Schur complement $J_{k,j,l}D_{k,j,l}^{-1}J_{k,j,l}^T$ is approximated by dropping each off-diagonal entry that falls below a relative tolerance $\tau_P$ and adding the magnitude of each dropped entry to the two corresponding diagonal entries, which preserves positive definiteness, and by adding $\delta_P$ to every diagonal entry, which amounts to replacing the zero lower right block of \eqref{eq.linsys} by $-\delta_P I$. On problems for which that approximation would be denser than the block system itself, the latter is factored directly instead. This preconditioner is symmetric but indefinite, so a GMRES-type method is employed.  Within \eqref{eq.linsys} we set $H_{k,j,l}$ to be the Hessian of the Lagrangian of \eqref{eq:sampled_slack} with respect to $x$ at the current primal-dual iterate, and we retain the primal-dual block $S_{k,j,l}Z_{k,j,l}$ rather than employ $\mu_{k,j}I$, subject to
\bequationNN
  s^{(i)}_{k,j,l} z^{(i)}_{k,j,l} \in [\underline\kappa_{sz}\mu_{k,j}, \overline\kappa_{sz}\mu_{k,j}]\ \ \text{for all}\ \ i \in [q],
\eequationNN
as in the default Ipopt implementation~\cite{wachter2006implementation}.

We apply the preconditioner described above on the right, within unrestarted flexible GMRES (FGMRES)~\cite{saad1993flexible}.  Writing \eqref{eq.linsys} as $\mathcal{A}w = h$, right preconditioning minimizes the unpreconditioned residual norm $\|h - \mathcal{A}w\|_2$, which is the norm bounded in \eqref{eq:kappa_des}, and FGMRES stores the preconditioned basis vectors, so that each iterate is available to the termination tests without an additional application of the preconditioner.  With left preconditioning, GMRES would instead minimize $\|\mathcal{P}^{-1}(h - \mathcal{A}w)\|_2$, where $\mathcal{P}$ denotes the preconditioner. Since FGMRES forms its iterates from the preconditioned basis vectors, which need not be well conditioned, their attainable accuracy can be limited.  We therefore compare the residual $\mathcal{R}_t := \|h - \mathcal{A}w_t\|_2$ of the iterate $w_t$, which the termination tests evaluate in any case, with the estimate $\widehat{\mathcal{R}}_t$ that FGMRES maintains in iteration $t$ and which equals it in exact arithmetic.  If $\mathcal{R}_t > \max\{\kappa_{\rm gap}\widehat{\mathcal{R}}_t, \epsilon_{\rm gap}\|h\|_2\}$, or if the Arnoldi process breaks down without an accepted iterate, the solve switches once to left-preconditioned GMRES~\cite{SaadSchu1986}, whose iterates are formed from an orthonormal basis, starting from the FGMRES iterate with the smallest $\mathcal{R}_t$.  In either phase, with $t$ counting the iterations of that phase, the solve is abandoned if $\mathcal{R}_t$ has failed to halve over the last $\max\{l^{\rm pd}_{\rm stall},\ \lfloor (t-1)/2 \rfloor\}$ iterations.  The iteration limit $l^{\rm pd}_{\max}$ applies to the two phases together.

To avoid an unnecessarily large number of inner iterations, Algorithm~\ref{alg.ipm} aims to produce reasonably accurate
solutions that satisfy the residual bound
\bequation\label{eq:kappa_des}
  \left\| \bbmatrix \rho_{k,j,l} \\ r_{k,j,l} \ebmatrix \right\|_2
  \leq \kappa_{\rm des}
  \left\| \bbmatrix g_{k,j,l} + J_{k,j,l}^Ty_{k,j,l} \\
                    -J_{k,j,l}v_{k,j,l} \ebmatrix \right\|_2
\eequation
for a given $\kappa_{\rm des} \in (0,1)$ before consideration of the termination tests; see~\cite{CurtScheWaec2010}.  This bound differs from the Residual Condition of \S\ref{sec.interior} in that it constrains both residual blocks and is relative to the current right-hand side alone.  However, if after a fixed number of iterations $\tilde l^{\rm pd}_{\max}$ this condition is not satisfied, then it is dropped and the algorithm requires only that one of Termination Tests~1, 2, or~3 is satisfied.  In evaluating these tests, the second inequality of the Tangential Step Condition is enforced with its right-hand side set to zero, as in~\cite{CurtHubeScheWaec2012}.  Since $\zeta\|v_{k,j,l}\|_2 \geq 0$, a step meeting this stronger requirement meets the condition as stated in \S\ref{sec.interior} for every $\zeta \in \R{}_{>0}$, so no value of $\zeta$ need be specified. The constants $\epsilon_2$ and $\epsilon_3$ in Table~\ref{tab:parameters} are those of Termination Tests~2 and~3 as stated in \S\ref{sec.interior}: $\epsilon_2$ determines when the dual-only step of TT2 may be considered, and $\epsilon_3$ is the fraction of the reduction in linearized infeasibility attained by the normal step that TT3 requires of the full step.  (In the numbering of~\cite{CurtHubeScheWaec2012}, the roles of these two constants are interchanged.)  In addition, for a step to be accepted through TT2, our implementation requires the Residual Condition to hold with a constant $\kappa_2 \in (0,\kappa]$ in place of $\kappa$, so that a dual-only step is accepted only when the dual residual has been reduced sufficiently.

Two formulations of \eqref{eq.linsys} are available in each iteration, and they differ as to whether a normal step is computed at all.  The \emph{direct} formulation replaces the lower right-hand-side block of \eqref{eq.linsys} by $-c_{k,j,l}$, so that $J_{k,j,l}d_{k,j,l} = -c_{k,j,l}$, which yields $d_{k,j,l}$ in a single solve and requires no normal step.  The \emph{decomposed} formulation, by contrast, first computes $v_{k,j,l}$ by Algorithm~\ref{alg:dogleg_normal}, then solves \eqref{eq.linsys} exactly as it is stated, thereby obtaining $d_{k,j,l} = v_{k,j,l} + u_{k,j,l}$.  The direct formulation is the less expensive of the two, since it avoids the normal-step solve entirely.  However, it offers no control over the tangential component when the iteration stalls: the Tangential Step Condition and the Model Reduction Condition are both stated in terms of $u_{k,j,l}$ and $v_{k,j,l}$ separately, and a short accepted step is the practical signal that these two quantities are in conflict.  We therefore switch between the two formulations, following~\cite[\S3.1]{CurtHubeScheWaec2012}.  Specifically, for a user-defined threshold $\kappa_{\rm sw} \in (0,1]$, an accepted step with $\alpha_{k,j,l} \geq \kappa_{\rm sw}$ selects the direct formulation for the subsequent iteration, whereas a shorter step selects the decomposed one, and each sampled problem begins with the direct formulation.  Each Krylov solve starts from the zero vector.

Two further situations merit attention.  First, the Hessian Modification Strategy is applied only once the same relative residual that~\eqref{eq:kappa_des} drives below $\kappa_{\mathrm{des}}$ has fallen to at most $\kappa_W$, since the curvature along an inaccurate iterate is not yet indicative of that of $W_{k,j,l}$ in the reduced space; see~\cite{CurtHubeScheWaec2012}.  Whenever the strategy calls for a modification of $W_{k,j,l}$, the current Krylov solve is discarded and a new one, with its own allowance of $l^{\rm pd}_{\max}$ iterations, is begun from the modified matrix.  Second, if an iterative solve ends without any of the three termination tests having been satisfied, whether by reaching $l^{\rm pd}_{\max}$ or through the safeguards above, then the trial step and all Hessian modifications accumulated along the way are discarded, and \eqref{eq.linsys} is instead solved by the inertia-correcting factorization stated as Algorithm~\ref{alg:inertia_correction}, which follows the default Ipopt algorithm~\cite{wachter2006implementation}, applied to the original $W_{k,j,l}$.  A capped direct solve is first retried once with the decomposed formulation and a freshly computed normal step.  The factorization applies only to the iteration in which it occurs; the next iteration begins with the iterative solver again.  Because the dual regularization perturbs the constraint block, a direction obtained in this way need not satisfy any of Termination Tests~1, 2, or~3 for \eqref{eq.linsys} as stated; it is accepted as a practical exception.

\begin{algorithm}[ht]
  \caption{Inertia Correction and Primal-Dual Solve}
  \label{alg:inertia_correction}
  \footnotesize
  \begin{algorithmic}[1]
    \State \textbf{Given} $W_{k,j,l}$, $J_{k,j,l}$, $g_{k,j,l}$, $y_{k,j,l}$, $v_{k,j,l}$, $\delta_w^0, \delta_c^0 \in \R{}_{>0}$
    \State Set $\delta_w \gets 0$, $\delta_c \gets 0$, $a_w \gets \delta_w^0$, and $a_c \gets \delta_c^0$
    \Loop
      \State Attempt to factor $K \gets \bbmatrix W_{k,j,l} + \delta_w I_{n+q} & J_{k,j,l}^T \\ J_{k,j,l} & -\delta_c I_{p+q} \ebmatrix$
      \If{the factorization fails}
        \State Set $\delta_w \gets a_w$, $\delta_c \gets a_c$, $a_w \gets 10a_w$, and $a_c \gets 10a_c$
      \Else
        \State Let $\varsigma_+$ and $\varsigma_-$ be the numbers of positive and negative eigenvalues of $K$, respectively
        \State \textbf{if} $\varsigma_+ = n+q$ \textbf{and} $\varsigma_- = p+q$ \textbf{then break}
        \State \textbf{if} $\varsigma_+ < n+q$ \textbf{then} $\delta_w \gets a_w$ and $a_w \gets 10a_w$
        \State \textbf{if} $\varsigma_- < p+q$ \textbf{then} $\delta_c \gets a_c$ and $a_c \gets 10a_c$
      \EndIf
    \EndLoop
    \State Solve $K\bbmatrix d_{k,j,l} \\ \delta_{y,k,j,l}\ebmatrix
      = -\bbmatrix g_{k,j,l} + J_{k,j,l}^Ty_{k,j,l} \\ -J_{k,j,l}v_{k,j,l}\ebmatrix$
      and \textbf{return} $(d_{k,j,l},\delta_{y,k,j,l})$
  \end{algorithmic}
\end{algorithm}

\subsection{Line search and iterate update}\label{sec.linesearch}

A few additional points should be mentioned about our implementation related to the line search and iterate update.  First, rather than perform a standard penalty parameter update and backtracking line search with sufficient decrease condition as stated in \eqref{eq:armijo}, our implementation employs a flexible penalty function as described in \cite{CurtNoce08}.  Specifically, the implementation maintains lower and upper merit parameter values $\pi^l_{k,j,l} \leq \pi^u_{k,j,l}$, initialized to $\pi^l_0$ and $\pi^u_0$, where the upper value plays the role of $\pi_{k,j,l}$ in Algorithm~\ref{alg.ipm}; in particular, it is the value increased through \eqref{eq.pi_trial}.  A step size $\alpha_{k,j,l} \in (0,\alpha^{\max}_{k,j,l}]$ is accepted if
\bequation\label{eq:flexible_armijo}
  \baligned
    &\ \phi_{\Scal_k}(x_{k,j,l} + \alpha_{k,j,l} \dtilde_{x,k,j,l}, s_{k,j,l} + \alpha_{k,j,l} \dtilde_{s,k,j,l}; \mu_{k,j}, \pi) \\
    \leq&\ \phi_{\Scal_k}(x_{k,j,l}, s_{k,j,l}; \mu_{k,j}, \pi) - \eta_2 \alpha_{k,j,l} \Delta m_{k,j,l}(d_{k,j,l}; \mu_{k,j}, \pi^u_{k,j,l})
  \ealigned
\eequation
for either $\pi = \pi^l_{k,j,l}$ or $\pi = \pi^u_{k,j,l}$.  If the accepted step satisfies \eqref{eq:flexible_armijo} for $\pi = \pi^u_{k,j,l}$ but not for $\pi = \pi^l_{k,j,l}$, then the lower value is increased toward the value at which the two sides of \eqref{eq:flexible_armijo} balance; namely, with $\varphi_{\Scal_k}(x,s;\mu) := f_{\Scal_k}(x) - \mu \sum_{i \in [q]} \ln(s^{(i)})$ and
\bequationNN
  \varrho_{k,j,l} := \frac{\varphi_{\Scal_k}(x_{k,j,l+1},s_{k,j,l+1};\mu_{k,j}) - \varphi_{\Scal_k}(x_{k,j,l},s_{k,j,l};\mu_{k,j})}{\|c_{k,j,l}\|_2 - \|c_{k,j,l+1}\|_2},
\eequationNN
the implementation sets
\bequationNN
  \pi^l_{k,j,l+1} \gets \min\left\{ \pi^u_{k,j,l},\ \pi^l_{k,j,l} + \max\{ \gamma_l (\varrho_{k,j,l} - \pi^l_{k,j,l}),\ \epsilon_l \} \right\}
\eequationNN
for user-defined $\gamma_l \in (0,1)$ and $\epsilon_l \in \R{}_{>0}$; otherwise, it sets $\pi^l_{k,j,l+1} \gets \pi^l_{k,j,l}$.  Moreover, the line search is capped at $\bar B$ trials.  If the line search does not terminate within this many trials, the final trial step is accepted rather than taking a null step.  Second, in place of the minimization \eqref{eq:dual_update_beta}, our implementation takes $\beta_{y,k,j,l}$ to be the least value in $[\alpha_{k,j,l},1]$ satisfying the inequality required of $y_{k,j,l+1}$ in \S\ref{sec.interior}, namely, \eqref{eq.mystery_inequality}.  Third, the tests for an infeasible stationary point in Algorithms~\ref{alg:ra} and~\ref{alg.ipm} are implemented: A point $(x,s)$ is declared an infeasible stationary point when $\|c_\Ecal(x)\|_\infty > \epsilon_{\rm infeas}$ or $\|\max\{c_\Ical(x),0\}\|_\infty > \epsilon_{\rm infeas}$ while the residual of \eqref{eq:infeasible_stationary} is at most $\epsilon_{k,1}$ in the $\ell_\infty$-norm.

\begin{table}[tp]
\centering
\footnotesize
\caption{Parameter values used for Algorithms~\ref{alg:ra}, \ref{alg.ipm},
  \ref{alg:dogleg_normal}, and~\ref{alg:inertia_correction}.}
\label{tab:parameters}
\setlength{\tabcolsep}{7pt}
\renewcommand{\arraystretch}{1.08}
\begin{tabular}{l l l}
\hline
\textbf{Symbol} & \textbf{Value} & \textbf{Description} \\
\hline

\multicolumn{3}{l}{\emph{Sampling (\S\ref{sec.impl_sampling})}} \\
\hline
$N$              & $1024$    & Full-sample size \\
$N_0$            & $128$     & Initial sample size \\
$\tilde\theta$   & $0.5$     & Sample-size accuracy parameter, \eqref{eq:N_k_var} \\
$\hat\beta$      & $5$       & Maximum sample-growth factor, \eqref{eq:N_k_rule} \\
\hline

\multicolumn{3}{l}{\emph{Warm-starting (\S\ref{sec.impl_warmstart})}} \\
\hline
$\mu_0$          & $10^{-1}$ & Initial barrier parameter, Alg.~\ref{alg:ra} \\
$\bar\rho$       & $1$       & Gradient-drift threshold, \eqref{eq:outer_restart_rel_drift} and subsequent discussion \\
$\underline\rho$ & $1$       & Drift denominator safeguard, \eqref{eq:outer_restart_rel_drift} \\
$c_\mu$          & $2$       & Barrier carry-forward factor, \eqref{eq:mu_hotstart} \\
$c_s$            & $10^{-2}$ & Slack residual-budget coefficient, \eqref{eq:slack_bump_pre} \\
\hline

\multicolumn{3}{l}{\emph{Normal step (\S\ref{sec.normal_impl_step})}} \\
\hline
$\omega$        & $100$     & Initial trust-region scaling, \eqref{eq:normal_step} \\
$\omega_{\max}$ & $10^{20}$ & Maximum trust-region scaling \\
$\epsilon_v$    & $1$       & Fraction of Cauchy decrease attained \\
$\delta_v$      & $10^{-8}$ & Initial normal-equation ridge \\
$\varsigma_v$   & $10$      & Normal-equation ridge growth factor \\
\hline

\multicolumn{3}{l}{\emph{Primal-dual step (\S\ref{sec.impl_step})}} \\
\hline
$\kappa$                  & $10^{-2}$ & Residual Condition constant \\
$\kappa_{\rm des}$        & $10^{-3}$ & Desired relative residual tolerance, \eqref{eq:kappa_des} \\
$\psi$                    & $0.1$     & Tangential Step Condition norm parameter \\
$\theta$                  & $10^{-12}\mu_{k,j}$ & Tangential Step Condition curvature parameter \\
$\sigma$                  & $0.09$    & Model Reduction Condition constant \\
$\epsilon_2$              & $10^{-8}$ & Eligibility tolerance of the dual-only test (TT2) \\
$\epsilon_3$              & $0.9$     & Feasibility-progress fraction of TT3 \\
$\kappa_2$                & $10^{-3}$ & Residual Condition constant for TT2, in place of $\kappa$ \\
$\kappa_W$                & $10^{-2}$ & Relative-residual threshold for Hessian modification, \S\ref{sec.impl_step} \\
$\delta_P$                & $10^{-10}$ & Preconditioner diagonal floor and Schur-complement shift \\
$\tau_P$                  & $10^{-3}$ & Preconditioner Schur drop tolerance \\
$\underline\kappa_{sz}$   & $10^{-10}$ & Lower safeguard on $s^{(i)}z^{(i)}/\mu_{k,j}$ \\
$\overline\kappa_{sz}$    & $10^{10}$  & Upper safeguard on $s^{(i)}z^{(i)}/\mu_{k,j}$ \\
$\kappa_{\rm sw}$         & $10^{-3}$ & Direct/decomposed switching threshold \\
$\delta_w^0 = \delta_c^0$ & $10^{-5}$ & Initial inertia corrections, Alg.~\ref{alg:inertia_correction} \\
$\tilde l^{\rm pd}_{\max}$ & $100$    & Krylov iterations before \eqref{eq:kappa_des} is dropped \\
$l^{\rm pd}_{\max}$       & $\lceil 1.1(n+p+2q) \rceil$ & Krylov iteration limit per linear system, both phases \\
$\kappa_{\rm gap}$        & $10$      & Residual-discrepancy factor for leaving FGMRES \\
$\epsilon_{\rm gap}$      & $10^{-14}$ & Absolute residual-discrepancy floor for leaving FGMRES \\
$l^{\rm pd}_{\rm stall}$  & $300$     & Minimum iterative stagnation window \\
\hline
\end{tabular}
\end{table}

\begin{table}[tp]
\centering
\footnotesize
\caption{Parameter values used for Algorithms~\ref{alg:ra}, \ref{alg.ipm},
  \ref{alg:dogleg_normal}, and~\ref{alg:inertia_correction} (continued).}
\label{tab:parameters2}
\setlength{\tabcolsep}{7pt}
\renewcommand{\arraystretch}{1.08}
\begin{tabular}{l l l}
\hline
\textbf{Symbol} & \textbf{Value} & \textbf{Description} \\
\hline

\multicolumn{3}{l}{\emph{Line search, barrier, and flexible penalty}} \\
\hline
$\eta_1$     & $\max\{0.99,\,1-\mu_{k,j}\}$ & Fraction-to-the-boundary parameter \\
$\eta_2$     & $10^{-8}$ & Sufficient-decrease constant, \eqref{eq:armijo} \\
$\tau$       & $0.1$     & Trial-penalty denominator factor, \eqref{eq.pi_trial} \\
$\pi^l_{0}$  & $10^{-8}$ & Initial lower flexible-penalty value \\
$\pi^u_{0}$  & $1$       & Initial upper flexible-penalty value \\
$\epsilon_\pi$ & $10^{-4}$ & Upper-penalty increment \\
$\gamma_l$   & $10^{-4}$ & Lower-penalty update fraction, \S\ref{sec.linesearch} \\
$\epsilon_l$ & $10^{-4}$ & Minimum lower-penalty increment, \S\ref{sec.linesearch} \\
$\epsilon_{\mu,a}$ & $10$ & Barrier-termination factor, $a \in [5]$, \eqref{eq:barrier_termination} \\
$\kappa_\mu$ & $0.2$     & Linear barrier-reduction factor, \eqref{eq:mu_update} \\
$\theta_\mu$ & $1.5$     & Superlinear barrier-reduction exponent, \eqref{eq:mu_update} \\
$\underline\mu$ & $10^{-9}$ & Barrier-parameter floor, \eqref{eq:mu_update} \\
\hline

\multicolumn{3}{l}{\emph{Iteration limits}} \\
\hline
$\bar K$           & $4$    & Maximum outer RA iterations \\
$\bar J$           & $8$    & Maximum barrier stages per sampled problem \\
$\bar L$           & $500$  & Maximum inner iterations per barrier subproblem \\
$\bar L_{\rm det}$ & $2000$ & Inner-iteration limit in deterministic mode \\
$\bar B$           & $50$   & Maximum line-search trials \\
\hline

\multicolumn{3}{l}{\emph{Termination tolerances, scaling, and test problems}} \\
\hline
$\gamma_{k,a}$ & $0.5$; $0$ once $N_k = N$ & Relative RA tolerance, $a \in [5]$, \eqref{eq:sample:termination} \\
$\epsilon_{k,a}$ & $10^{-5}$ & Absolute RA tolerance, $a \in [5]$, \eqref{eq:sample:termination} \\
$\epsilon_{\rm full}$ & $10^{-4}$ & Full-sample termination tolerance \\
$\epsilon_{\rm infeas}$ & $10^{-4}$ & Infeasibility threshold for infeasible stationarity, \S\ref{sec.linesearch} \\
$\Delta_\xi$ & $10^{-1}$ & Noise half-width, \eqref{eq:augmented_objective} \\
$g_{\max}$   & $100$     & Gradient cap for problem scaling \\
$s^{\rm scale}_{\min}$ & $10^{-8}$ & Minimum scaling factor, \S\ref{sec.prescaling} \\
\hline

\end{tabular}
\end{table}

\section{Numerical Experiments}\label{sec.numerical}

Our aim in this section is to demonstrate the practical performance of our implementation of RA-IPM described in \S\ref{sec.implementation} through three sets of numerical experiments.  The first tests the effect of the warm-starting strategy between sampled problems, the second demonstrates the improved performance of the RA-IPM algorithm due to adaptive sampling through retrospective approximation as compared with solving the full-sample problem directly, and the third shows the improvement in solution time that is possible due to the use of an interior-point method rather than an SQP method, specifically when solving larger-scale problems.  Since the benefits of inexact search direction computations are only relevant for such larger-scale problems, in the first two experiments we employ a factorization-based approach for computing the search directions, and only employ the iterative linear system solves with inexact computations in the third experiment.

All of the results in this section were obtained on a MacBook Air with an Apple M2 processor, eight CPU cores, and 8~GB of memory running macOS~26.5.1.  Our implementation is written in Python~3.11.11 with NumPy~2.4.2~\cite{harris2020array} and SciPy~1.17.1~\cite{virtanen2020scipy}, the latter providing the sparse LU factorization used for the preconditioner of \S\ref{sec.impl_step}.  Our test problems, described below, are accessed through PyCUTEst~1.8.0~\cite{pycutest} with CUTEst~2.6.0~\cite{cutest}.  The normal equations of Algorithm~\ref{alg:dogleg_normal} and the symmetric indefinite systems of Algorithm~\ref{alg:inertia_correction} are factored with MUMPS~5.8.2~\cite{amestoy2001mumps} through python-mumps~0.0.6.  Our experiments also involve an implementation of RA-SQP~\cite{berahas2025retrospective}, which uses Gurobi~13.0.2~\cite{gurobi} through gurobipy~13.0.2.  We also show the results of runs using Ipopt as the subproblem solver; specifically, we employed Ipopt~3.14.19~\cite{wachter2006implementation}, accessed through CyIpopt~1.6.1 with MUMPS as its linear solver.

The parameters of Algorithms~\ref{alg:ra} and~\ref{alg.ipm} were set as stated in Tables~\ref{tab:parameters}--\ref{tab:parameters2}, except where an experiment states otherwise.  Ipopt departs from its defaults in that its own scaling is disabled in favor of the scaling we implemented as described in the beginning of~\S\ref{sec.implementation}.  We also set its iteration limit to match that of RA-IPM and modified its stopping test through its iteration callback so that it terminates on the conditions of \eqref{eq:sample:termination}. The remaining settings were left at Ipopt defaults.

We generated our test set from the CUTEst collection~\cite{cutest} through the PyCUTEst interface~\cite{pycutest}, proceeding as follows.  Firstly, we retained the problems whose constraints are classified as quadratic or more general, and removed all others.  Secondly, we expressed each such problem in the form \eqref{eq:full} and removed those possessing no inequality constraint that is nonlinear since RA-IPM is primarily designed for problems with at least one inequality constraint that is nonlinear.  Thirdly, we removed those problems with $n+p+q > 7000$.  Finally, we removed those problems whose objective, constraint, or derivative values are not finite at the prescribed starting point, namely LHAIFAM, S365, and S365MOD.  Our test set was the remainder of these problems, which resulted in $169$ total problems.

The experiments of \S\ref{sec.exp_hotstart} and \S\ref{sec.exp_sampling} measure component evaluations, where each evaluation of $\nabla f_{\Scal_k}$ during these experiments costs $N_k$ gradient component evaluations. The experiment in \S\ref{sec.exp_large}, on the other hand, measures CPU time.  For the latter, the dominant cost is that of solving \eqref{eq.linsys} in each iteration. As an estimator for that cost, we considered the number of nonzeros (nnz) of the matrix.  We accordingly measured the size of a problem $P$ by the number of nonzeros in
\bequationNN
  K_0(P) := \bbmatrix W_0 & J_0^T \\ J_0 & 0 \ebmatrix
  \in \R{(n+p+2q) \times (n+p+2q)},
\eequationNN
where $W_0$ and $J_0$ are the matrices of \S\ref{sec.interior} evaluated at the initial primal-dual iterate of the scaled problem with barrier parameter $\mu_0$. Let us refer to this count as $\mathrm{nnz}(K_0(P))$.  We defined a \emph{small} set $\Pcal_{\rm s}$ of problems, employed in the first two experiments, consisting of the $126$ (out of 169) problems with $\mathrm{nnz}(K_0(P)) \leq 2000$, whereas we defined a \emph{large} set $\Pcal_{\rm l}$ of problems consisting of the remaining $43$ problems, which we employed in the third experiment.

Following \cite{berahas2025retrospective}, we modified the default CUTEst problems with an additional objective term in order to allow for sampling through retrospective approximation.  Letting~$f_{\rm orig}$ denote the original CUTEst objective function, $x_{\rm init}$ denote the prescribed starting point, and $e_n \in \R{n}$ denote a vector of ones, each objective is modified to be
\bequation\label{eq:augmented_objective}
  f(x) = f_{\rm orig}(x) + \frac1N \sum_{i=1}^N \xi_i \|x - x_{\rm init} - e_n\|_2^2,
\eequation
where $\xi_i$ for each $i \in [N]$ is drawn from a uniform distribution over $[-\Delta_\xi,\Delta_\xi]$.  Each sampled objective $f_\Scal$ is formed by choosing uniformly, at random, from $\{\xi_i\}_{i\in[N]}$ without replacement. However, since such a selection might lead to $\frac1{|\Scal|} \sum_{i\in\Scal} \xi_i$ being negative, a selection is thrown out and re-drawn to ensure that $\frac1{|\Scal|} \sum_{i\in\Scal} \xi_i > 0$. The same goes for the full-sample problem to ensure that $\frac1N \sum_{i=1}^N \xi_i > 0$.  Each method that was employed to solve a given problem received the same seed, which ensured that, when comparing different methods for solving the same problem, the full sample $\Scal_{\rm full}$ was the same across methods.

Finally, let us state how we measure performance.  A run is successful when, at some iterate, the residuals on the left-hand side of \eqref{eq:sample:termination}, evaluated with the full-sample gradient, satisfy $\max_{a \in [5]} \|r_a\|_\infty \leq \epsilon_{\rm full} := 10^{-4}$.  All five residuals are employed for runs employing RA-IPM and Ipopt.  RA-SQP, on the other hand, maintains no multipliers for the bounds on the slack variables, so for it the criterion is applied with residuals $r_1$, $r_3$, $r_4$, and the norm of $S y_\Ical$.  Runs that exhaust an iteration limit, that stop at an infeasible stationary point, or that fail numerically are counted as unsuccessful.  We report the performance using Dolan--Mor\'e performance profiles~\cite{dolan2002benchmarking}. Thus, for each problem the measure of a method is divided by the smallest measure attained on that problem by any method under comparison, where an unsuccessful run is assigned an infinite ratio. Each method is run once per problem and configuration.

\subsection{Warm-starting between sampled problems}\label{sec.exp_hotstart}

Our first set of experiments examines the warm-starting strategy of \S\ref{sec.impl_warmstart}.  We compare the \emph{warm} configuration, which applies \eqref{eq:mu_hotstart} and \eqref{eq:slack_bump}, against a \emph{cold} configuration, which in every outer iteration sets $\mu_{k,0} = \mu_0$, restores the full budget of $\bar J$ barrier stages, and sets $\iota_{s,k} = 0$ and $\iota_{z,k} = 0$.  In both configurations the primal-dual iterate returned by each inner solve is carried forward, as line~\ref{step.slack_perturbation} of Algorithm~\ref{alg:ra} prescribes, so the two runs differ only in the barrier parameter, the barrier-stage budget, and the slack variables with which each inner solve begins.  Both are run over $\Pcal_{\rm s}$ with the same parameter values.  We report profiles of function component evaluations, gradient component evaluations, and inner iterations.

\begin{figure}[htbp]
  \centering
  \includegraphics[width=0.32\textwidth]{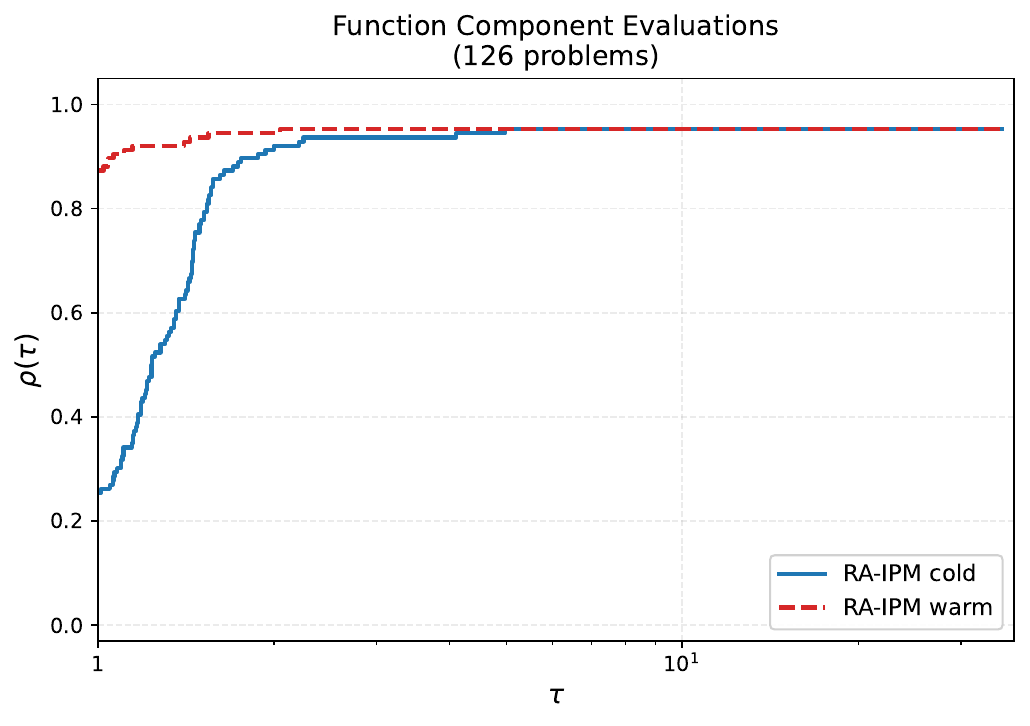}
  \includegraphics[width=0.32\textwidth]{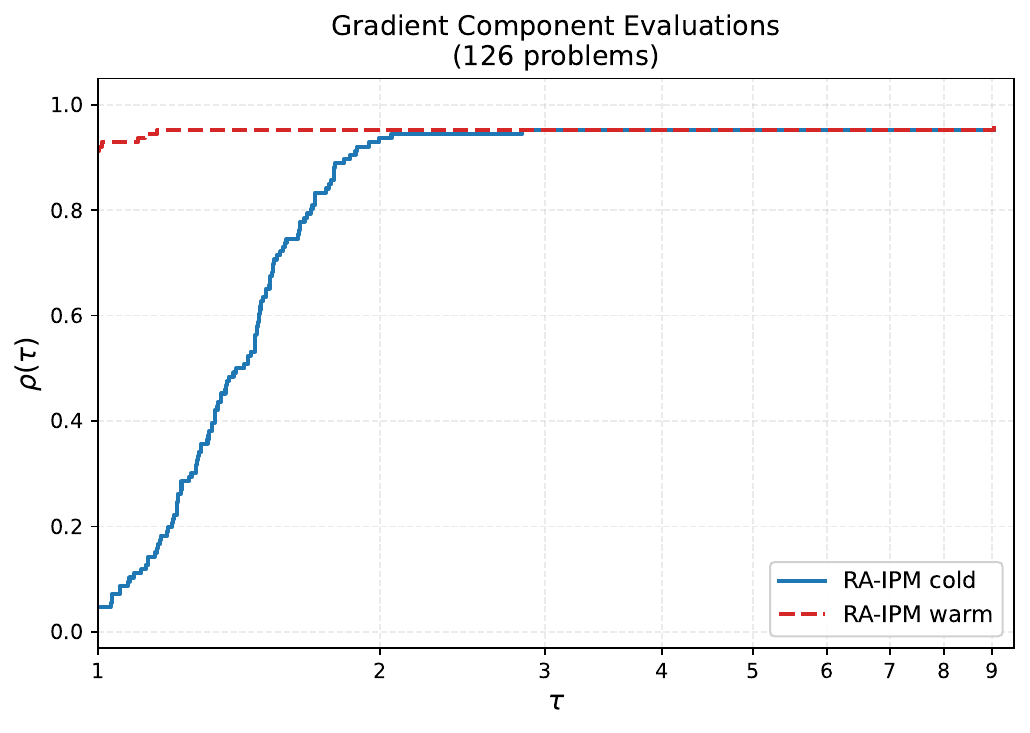}
  \includegraphics[width=0.32\textwidth]{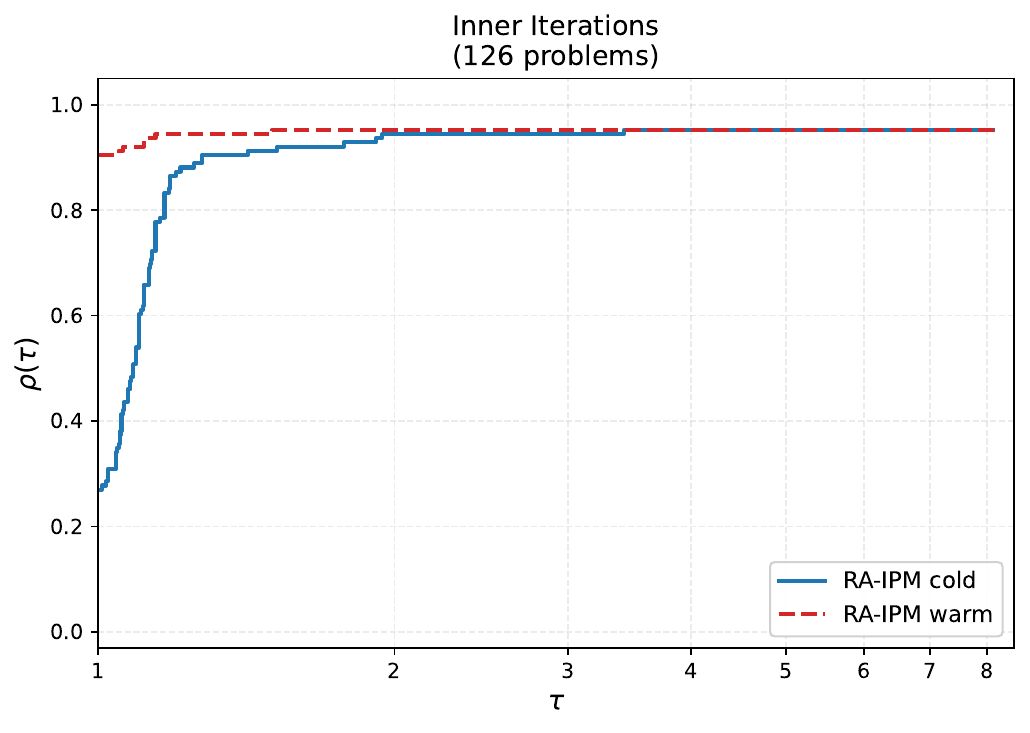}
  \caption{Performance profiles for the warm and cold configurations of RA-IPM
  on the $126$ problems of $\Pcal_{\rm s}$: function component evaluations,
  gradient component evaluations, and inner iterations.}
  \label{fig:contribution1}
\end{figure}

Figure~\ref{fig:contribution1} provides the results.  The warm configuration solves $121$ of the $126$ problems and the cold configuration $120$, all of which are also solved by the warm configuration.  For all three measures the warm profile starts above the cold one at a performance ratio of one, but the size of the effect differs.  It is clearest for gradient component evaluations: on $115$ of the $120$ problems solved by both configurations, the warm configuration uses fewer of them.  For function component evaluations and inner iterations, the cold profiles catch up within a small factor.

This pattern reflects the way in which the cold configuration restarts.  Its first barrier subproblem in each outer iteration uses $\mu_0$, whereas the iterate carried forward is nearly stationary for the previous sampled problem.  The barrier termination test \eqref{eq:barrier_termination}, with $\epsilon_{\mu,a} = 10$, holds at such a point for every sufficiently large value of the barrier parameter, so several barrier stages terminate at their entry test without taking a step.  Each such stage costs a sampled-gradient evaluation for its entry test but records no inner iteration.  Warm-starting therefore saves gradient evaluations mainly by sparing these stages, and saves comparatively few inner iterations.  The slack perturbation \eqref{eq:slack_bump} is small by construction, with $\|\iota_{s,k}\|_\infty \leq c_s \epsilon_{k,4} = 10^{-7}$ for the values in Tables~\ref{tab:parameters}--\ref{tab:parameters2}, and the experiment does not separate its effect from that of the barrier parameter.

\subsection{Adaptive sampling against deterministic solves}\label{sec.exp_sampling}

The paradigm of RA rests upon the premise that the early outer iterations do not require accurate objective information.  Should this premise hold, a method that begins with $N_0 = 128$ realizations and enlarges the sample only as the residual becomes small ought to reach the target accuracy with fewer component evaluations than one that employs all $N$ realizations from the start.  The premise is not automatic, however, since enlarging the sample alters the problem being solved, so that part of the work performed over $\Scal_{k-1}$ is discarded.  Our second set of experiments tests this premise, and does so for three inner solvers, so that any advantage observed may be attributed to the outer loop rather than to the subproblem solver.

We compare the \emph{adaptive} mode of \S\ref{sec.impl_sampling} against a \emph{deterministic} mode in which the outer loop receives $\Scal_{\rm full}$ in its first iteration.  In deterministic mode, one has $N_0 = N$, so that by the rule of \S\ref{sec.impl_sampling} the relative tolerances $\gamma_{k,a}$ vanish and only the absolute conditions in \eqref{eq:sample:termination} remain. In adaptive mode, the relative tolerances are likewise set to zero once the sample reaches $\Scal_{\rm full}$.  In both modes, each method also evaluates the residuals with the full-sample gradient in every inner iteration and stops as soon as the success criterion stated at the beginning of this section holds.  These full-sample evaluations serve only to detect success and are not counted. That said, all other evaluations are counted, including those used to estimate the variance in \eqref{eq:Z_V_estimates}.  The measures reported below are therefore the costs of reaching the target accuracy $\epsilon_{\rm full}$, which in adaptive mode may be reached before the sample has grown to $\Scal_{\rm full}$.  We note that this stopping rule relies on the full-sample gradient, which would not be available in the setting of \S\ref{sec.problem}. In that setting, the method would instead continue until $N_k = N$ and the absolute tolerances $\{\epsilon_{k,a}\}$ are met.  We employ the rule, identically in both modes, so that the two modes are compared on the cost of attaining the same accuracy for the same instance of \eqref{eq:full}. Measured instead at each mode's own termination, the adaptive mode would also be charged for the full-sample solve that follows.

The three inner solvers are RA-IPM, RA-SQP~\cite{berahas2025retrospective}, and Ipopt~\cite{wachter2006implementation}.  Ipopt is driven through its iteration callback and terminated by the conditions \eqref{eq:sample:termination}, as RA-IPM is; at the start of each sampled problem it is started from the current primal iterate with initial barrier parameter $\mu_0$ and its default initialization of the multipliers.  RA-SQP terminates each sampled solve by its step-norm test~\cite{berahas2025retrospective}.  In each case the inner solver is applied to the sampled problem within Algorithm~\ref{alg:ra} and receives the same $\Scal_{\rm full}$ in both modes, so that each pair of profiles constitutes a comparison within a single inner solver.  In particular, the deterministic runs are not stock Ipopt or stock RA-SQP applied to a different problem, and we draw no conclusion across the three pairs.  The iteration limits are $\bar L = 500$ in adaptive mode and $\bar L_{\rm det} = 2000$ in deterministic mode for both RA-IPM and RA-SQP, and $3000$ for Ipopt in both modes.  The units differ: for RA-IPM the limit applies to each barrier subproblem, whereas for RA-SQP and Ipopt it applies to each sampled solve.

We report function and gradient component evaluations, as well as total inner iterations.  A component evaluation is an evaluation of $F(\cdot,\xi)$, or of its gradient, for one realization~$\xi$, so that an evaluation of $f_{\Scal_k}$ or of $\nabla f_{\Scal_k}$ is charged $N_k$ components.  One expects adaptive sampling to prevail on the first two measures, since its early evaluations involve only a fraction of the full sample, and to fare worse on the third, since it may solve several problems rather than one.  For RA-IPM and Ipopt the third measure counts inner solver iterations, whereas for RA-SQP it counts the iterations consumed by Gurobi within the subproblem solves.

\begin{figure}[htbp]
  \centering
  \includegraphics[width=0.32\textwidth]{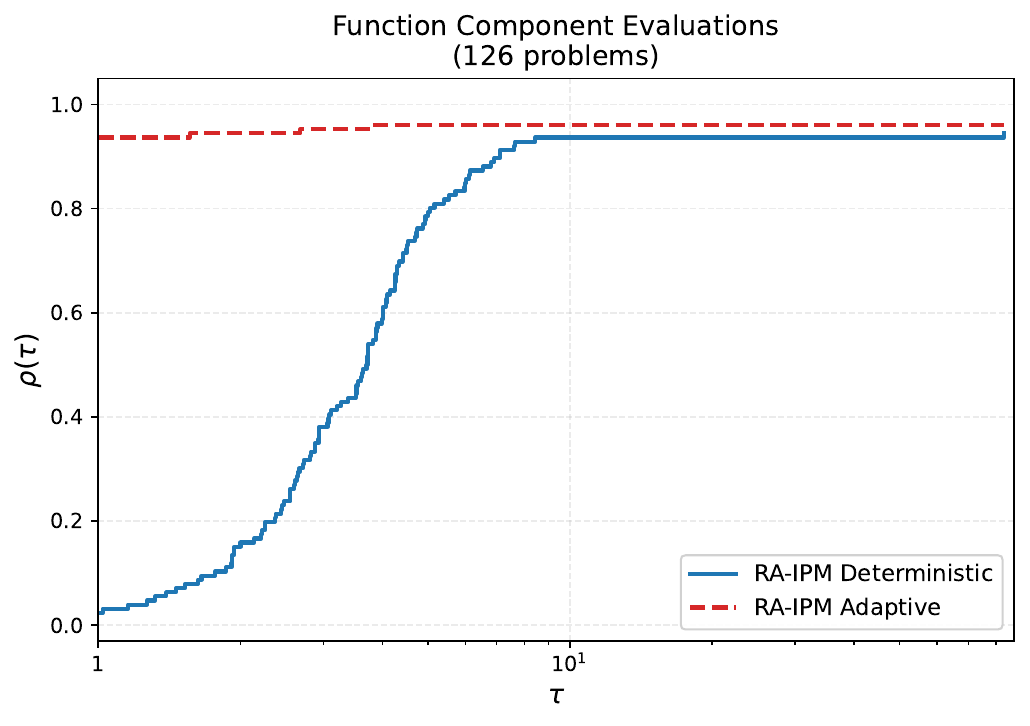}
  \includegraphics[width=0.32\textwidth]{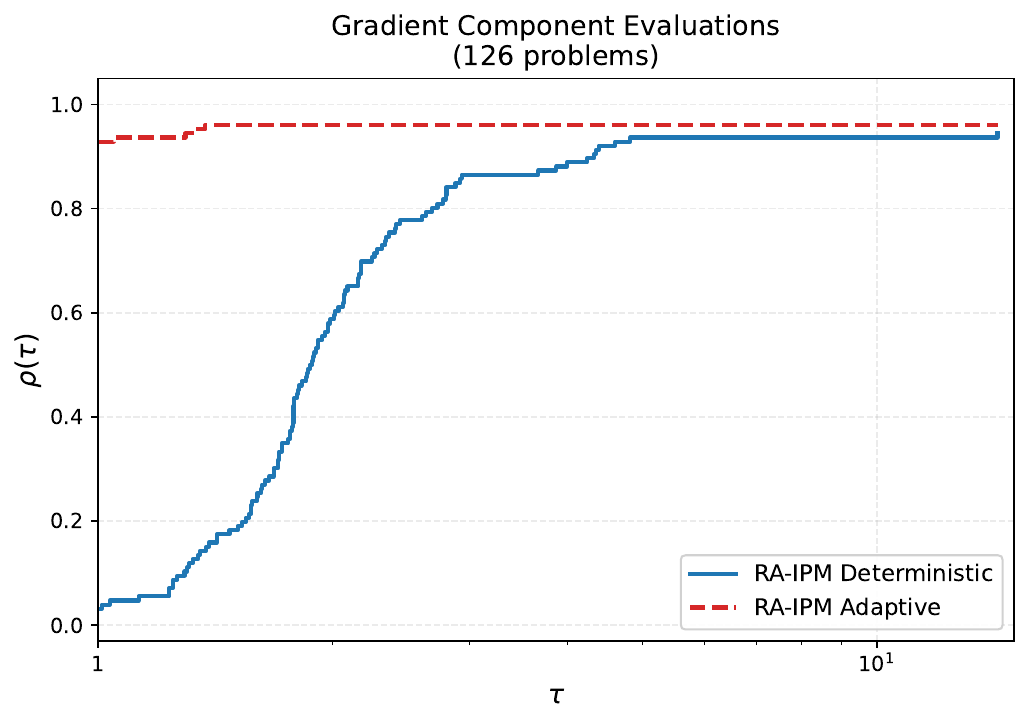}
  \includegraphics[width=0.32\textwidth]{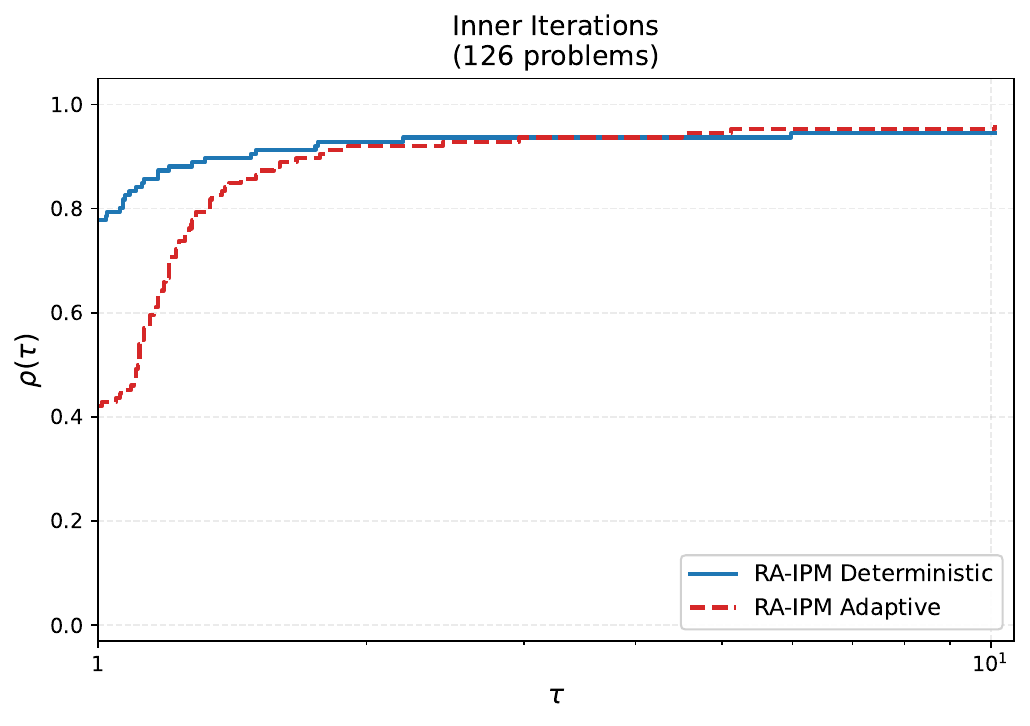}
  \caption{Performance profiles for the deterministic and adaptive modes with
  RA-IPM on the $126$ problems of $\Pcal_{\rm s}$: function component evaluations, gradient component
  evaluations, and inner iterations.}
  \label{fig:contribution2}
\end{figure}

\begin{figure}[htbp]
  \centering
  \includegraphics[width=0.32\textwidth]{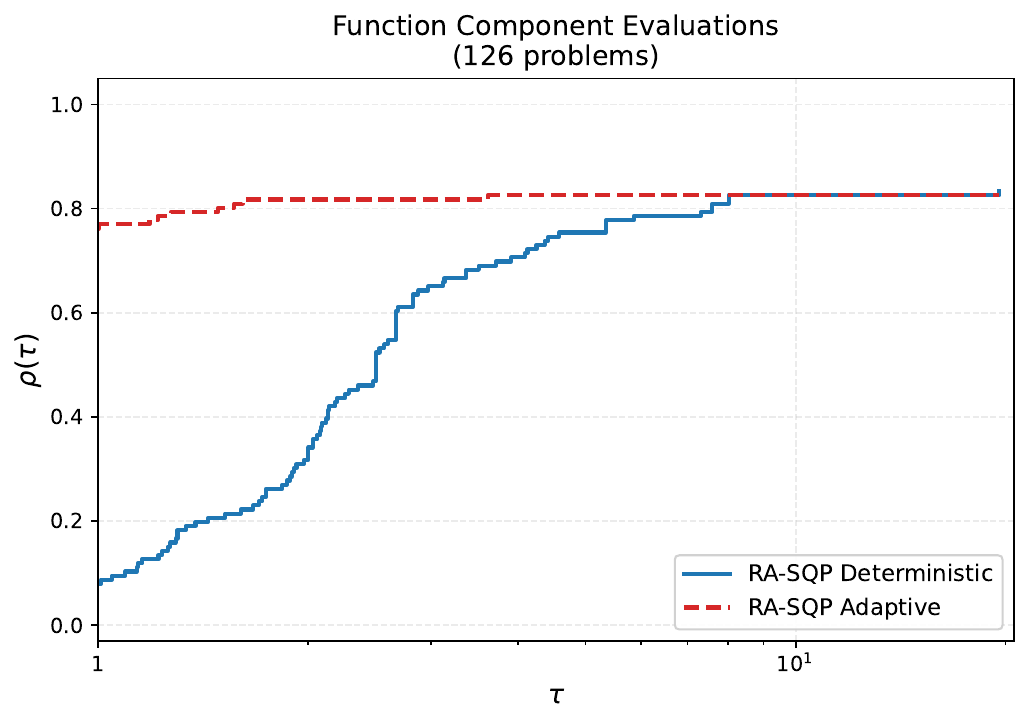}
  \includegraphics[width=0.32\textwidth]{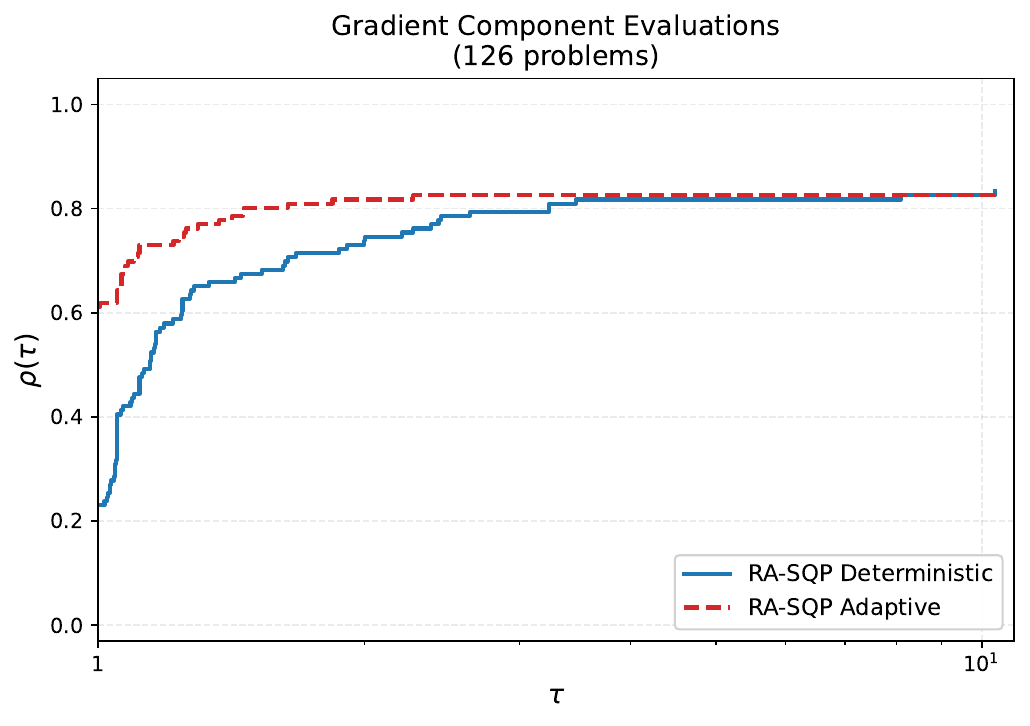}
  \includegraphics[width=0.32\textwidth]{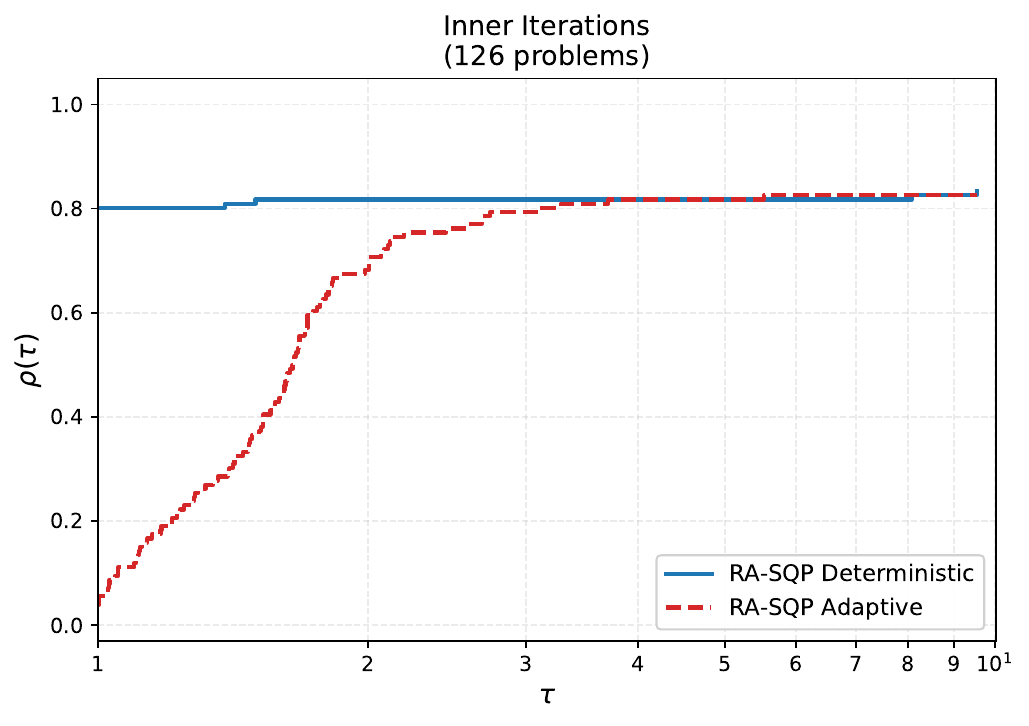}
  \caption{Performance profiles for the deterministic and adaptive modes with
  RA-SQP on the $126$ problems of $\Pcal_{\rm s}$: function component evaluations, gradient component
  evaluations, and inner iterations, here counted as Gurobi subproblem iterations.}
  \label{fig:contribution2_sqp}
\end{figure}

\begin{figure}[htbp]
  \centering
  \includegraphics[width=0.32\textwidth]{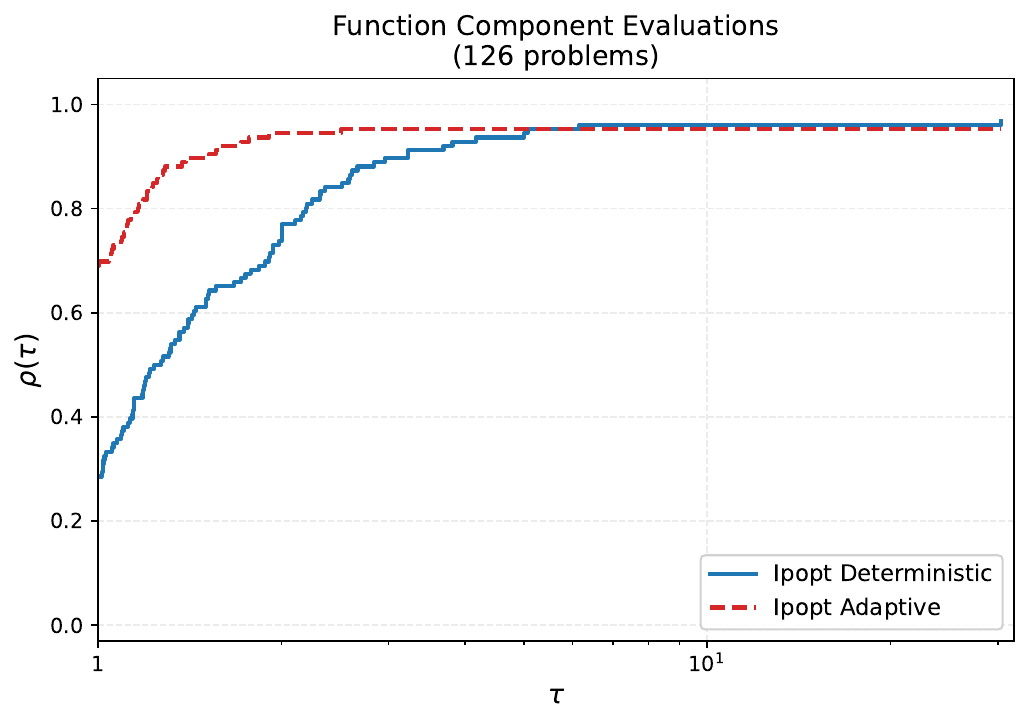}
  \includegraphics[width=0.32\textwidth]{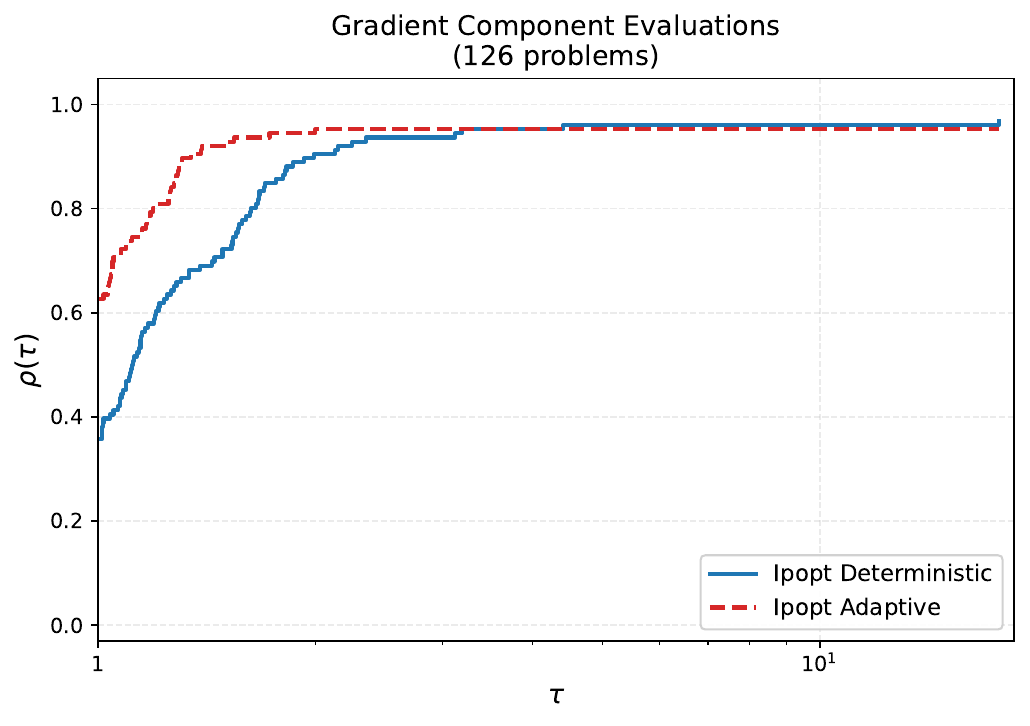}
  \includegraphics[width=0.32\textwidth]{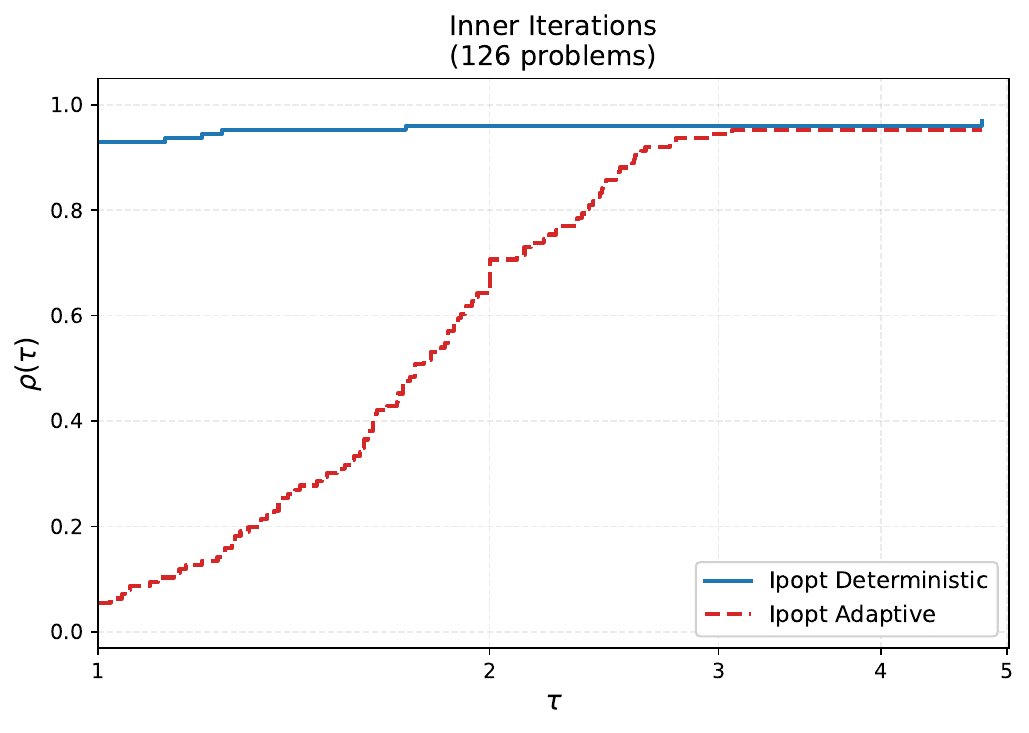}
  \caption{Performance profiles for the deterministic and adaptive modes with
  Ipopt on the $126$ problems of $\Pcal_{\rm s}$: function component evaluations, gradient component
  evaluations, and inner iterations.}
  \label{fig:contribution2_ipopt}
\end{figure}

Figures~\ref{fig:contribution2}--\ref{fig:contribution2_ipopt} provide the results.  For each of the three inner solvers, the adaptive profiles for component evaluations start well above the deterministic ones at a performance ratio of one, whereas the deterministic profiles approach their final levels only at larger ratios.  With RA-IPM the adaptive profile is flat from a ratio of one, so the adaptive mode attains the smaller count on essentially every problem it solves; with RA-SQP and Ipopt the advantage is less uniform.  The inner-iteration profiles, by contrast, favor the deterministic mode for each pair, so the savings come from the cost of each evaluation rather than from the number of iterations, which is the pattern that the RA paradigm predicts.  The two modes are similarly reliable, the adaptive mode solving two problems more than the deterministic mode with RA-IPM, one fewer with RA-SQP, and two fewer with Ipopt.

\subsection{Problems with larger KKT systems}\label{sec.exp_large}

Over our set of larger problems, the matrix in~\eqref{eq.linsys} has more than $2000$ and at most $4{,}007{,}000$ nonzeros at the starting point.  These problems are not of small size, so this experiment concerns the regime in which the linear algebra of a step begins to dominate the cost of a run.  Our set of experiments here compares RA-IPM against RA-SQP~\cite{berahas2025retrospective} on these problems.  For deterministic problems, the advantage of interior-point methods over SQP methods as problems grow is generally recognized.  The question here is whether this advantage persists when the solver must be restarted for each new sampled problem.

The two methods compute a step quite differently.  An RA-IPM iteration requires an approximate solution of one linear system of the form \eqref{eq.linsys}, obtained by the Krylov procedure of \S\ref{sec.impl_step}, with a fallback to the inertia-correcting factorization of Algorithm~\ref{alg:inertia_correction} when no candidate direction satisfies a termination test.  Such an iteration still involves sparse factorizations, namely of the approximate Schur complement used by the preconditioner and, in the decomposed formulation, of the matrix of the normal equations in Algorithm~\ref{alg:dogleg_normal}, but the matrix in \eqref{eq.linsys} itself is factored only when the fallback is invoked.  An RA-SQP iteration, by contrast, solves a dual linear program, a feasibility linear program, and a quadratic program with Gurobi, and computes a smallest-eigenvalue estimate for the QP Hessian matrix.  For our experiments here, both methods are run with the same $\Scal_{\rm full}$, the same seed, and the same sampling rule, although the samples they draw need not agree since the sample sizes depend on the iterates through the rule of \S\ref{sec.impl_sampling}.

Over $\Pcal_{\rm l}$ the time consumed by a run is dominated by the linear algebra of the step computation rather than by the evaluation of $f_{\Scal_k}$ and $\nabla f_{\Scal_k}$, so component evaluation counts no longer discriminates between the methods.  Iteration counts are not comparable either, since the two methods perform very different work per iteration.  We therefore compare CPU time, measured as the processor time of each run summed over all of its threads.  NumPy, SciPy, and MUMPS were restricted to a single thread, Gurobi was run with its default settings, and no time limit was imposed.  CPU time compares implementations as well as algorithms.  RA-SQP solves its subproblems with Gurobi, a mature commercial code, but each of its iterations also incurs the cost of building those subproblems and of the eigenvalue estimate. RA-IPM, on the other hand, is research code and uses publicly available linear algebra routines.  Nonetheless, the comparison should be read as indicating the relative cost of the two approaches as implemented here, rather than as a benchmark of solvers.

\begin{figure}[htbp]
  \centering
  \includegraphics[width=0.7\textwidth]{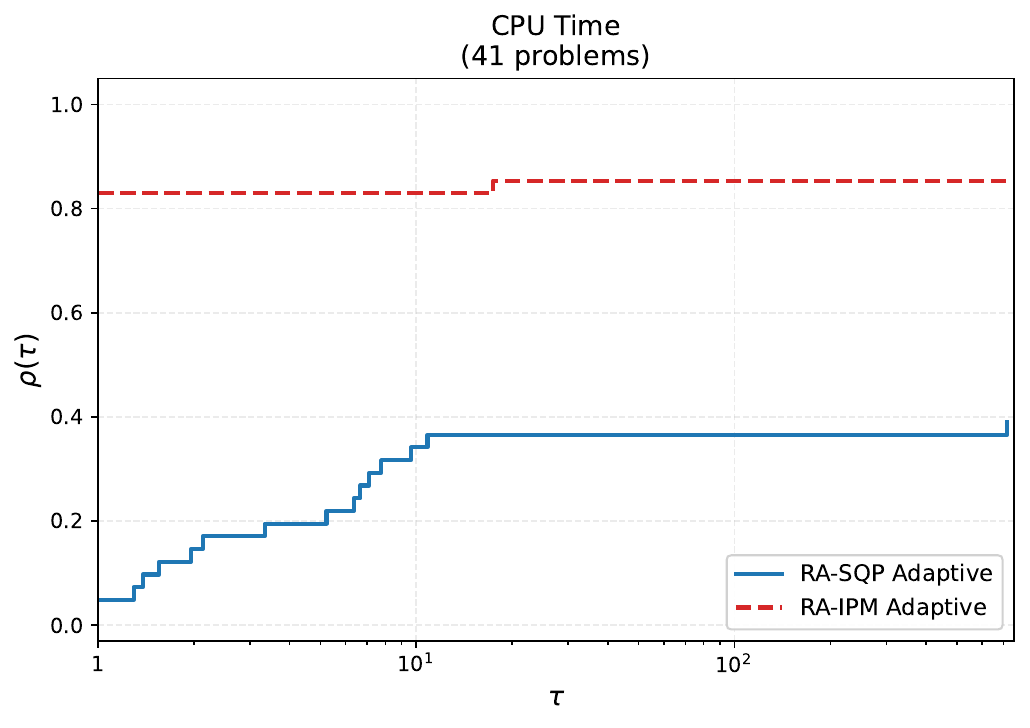}
  \caption{Performance profile of CPU time for RA-SQP and RA-IPM, both in
  adaptive-sampling mode, on $41$ of the $43$ problems of $\Pcal_{\rm l}$; see
  the text for the remaining two.}
  \label{fig:contribution3}
\end{figure}

Figure~\ref{fig:contribution3} reports results for $41$ of the $43$ problems of $\Pcal_{\rm l}$. (On the remaining two, the runs of both methods exhausted the available memory, so we omit them.)  RA-IPM solves $35$ of the $41$ problems while RA-SQP only solves $16$.  Fifteen problems are solved by both methods, $20$ only by RA-IPM, and one only by RA-SQP.  RA-IPM attains the smallest CPU time on $34$ of the $41$ problems and RA-SQP on $2$.  This ordering reflects an advantage in time and not only in reliability, since on the $15$ problems solved by both methods, RA-IPM is the faster on $14$.  On this test set and with our implementations, RA-IPM was thus both more reliable and faster than RA-SQP.

\section{Conclusion}\label{sec.conclusion}

We have proposed RA-IPM, an algorithm for minimizing an expectation subject to deterministic nonlinear equality and inequality constraints.  Its outer loop follows the retrospective-approximation paradigm, solving a sequence of sample-average problems with increasing sample sizes, and each sampled problem is solved approximately by an interior-point method that allows inexact step computations.  The inner method stops on the same conditions under which the outer loop is analyzed, namely a relative reduction, up to an absolute tolerance, of each block of the KKT residual.

Under our assumptions, if the algorithm generates an infinite sequence of iterates, the sample sizes grow without bound, and the tolerances vanish, then the expected KKT residuals at the outer iterates converge to zero.  These assumptions include a bound on the sampling error of the gradient at the points returned by the inner solves.  The analysis is asymptotic and does not bound the total work of the method.

In experiments on CUTEst problems whose objectives were augmented with a sampled term, carrying the barrier parameter and the slack variables from one sampled problem to the next reduced the number of gradient component evaluations relative to restarting the barrier parameter, mainly by sparing barrier stages.  For each of three inner solvers, solving a sequence of sampled problems of increasing size reduced the number of component evaluations needed to reach the target accuracy on the full-sample problem, at the price of more inner iterations.  On problems with larger KKT systems, RA-IPM solved more problems than RA-SQP and was faster on all but one of the problems that both solved.

\bibliographystyle{plain}
\bibliography{references}

\end{document}